\documentclass[11pt]{article}          
\usepackage{graphicx}
\usepackage{hyperref}
\usepackage{amsmath}
\usepackage{fullpage}
\usepackage{amssymb}
\usepackage{amsfonts}
\usepackage{latexsym}
\usepackage{amsthm}

\usepackage{lipsum}
\usepackage{amsfonts}
\usepackage{graphicx}
\usepackage{epstopdf}
\usepackage{algorithmicx}
\usepackage{algorithm}
\newtheorem{example}{Example}
\usepackage{pbox}
\ifpdf%
  \DeclareGraphicsExtensions{.eps,.pdf,.png,.jpg}
\else
  \DeclareGraphicsExtensions{.eps}
\fi
\usepackage{amsopn}

\usepackage{booktabs}

\DeclareMathOperator*{\argmin}{argmin}       
\DeclareMathOperator*{\argmax}{argmax}

\usepackage{float}
\usepackage{breqn}
\usepackage{multirow}
 \usepackage{placeins}

\newcommand{\hf}{\frac12}

\newcommand{\bfv}{ {\bf{v}}}

\newcommand{\bbR}{\mathbb{R}}

\newcommand{\bfA}{{\bf A}}
\newcommand{\bfD}{{\bf D}}

\newcommand{\bfJ}{{\bf J}}
\newcommand{\bfL}{{\bf L}}

\newcommand{\bfV}{{\bf V}}
\newcommand{\bfI}{{\bf I}}

\newcommand{\bfW}{{\bf W}}
\newcommand{\bfH}{{\bf H}}

\newcommand{\bfd}{{\bf d}}

\newcommand{\bfs}{{\bf s}}

\newcommand{\bfr}{{\bf r}}

\newcommand{\bfy}{{\bf  y}}
\newcommand{\bfx}{{\bf  x}}

\newcommand{\bfu}{{\bf u}}
\newcommand{\bfz}{{\bf z}}

\newcommand{\bfp}{{\bf p}}

\newcommand{\bfGamma}{\boldsymbol{\Gamma}}

\newcommand{\bfepsilon}{\boldsymbol \epsilon}

\graphicspath{{figures/}}

\date{\today}

\newcommand*\samethanks[1][\value{footnote}]{\footnotemark[#1]}
\begin{document}

\title{An Uncertainty-Weighted Asynchronous ADMM Method for Parallel PDE Parameter Estimation} 



\author{%
Samy Wu Fung\thanks{Department of Mathematics and Computer Science, Emory University,  Atlanta, GA, USA. \texttt{\{samy.wu,lruthotto\}@emory.edu}} \and Lars Ruthotto\samethanks[1]}
\maketitle
\begin{abstract}
We consider a global variable consensus ADMM algorithm for solving large-scale PDE parameter estimation problems asynchronously and in parallel. To this end, we partition the data and distribute the resulting subproblems among the available workers. Since each subproblem can be associated with different forward models and right-hand-sides, this provides ample options for tailoring the method to different applications including multi-source and multi-physics PDE parameter estimation problems. We also consider an asynchronous variant of consensus ADMM to reduce communication and latency.

Our key contribution is a novel weighting scheme that empirically increases the progress made in early iterations of the consensus ADMM scheme and is attractive when using a large number of subproblems. This makes consensus ADMM competitive for solving PDE parameter estimation, which incurs immense costs per iteration. The weights in our scheme are related to the uncertainty associated with the solutions of each subproblem. We exemplarily show that the weighting scheme combined with the asynchronous implementation improves the time-to-solution for a 3D single-physics and multiphysics PDE parameter estimation problems.

\vspace{3mm}
\end{abstract}

\section{Introduction}\label{sec:intro}
Recent technological advances have allowed us to collect data at massive scales with relative ease. This trend, often referred to as 'Big Data', has given rise to notoriously challenging high-dimensional parameter estimation problems. A common example is the computation of the {maximum a posteriori} (MAP) estimate \cite{calvetti2006large,stuart2010inverse} of large-scale Bayesian inverse problems. In this case, the parameter estimation problem is solved iteratively using gradient-based optimization. This requires numerous simulations that may involve different physical models and commonly scale to millions of variables \cite{biegler2003large}, leading to very high costs in both CPU time and memory. As a result, parallel and distributed iterative optimization techniques have become highly desirable, if not necessary, for solving these types of problems. 

In this paper, we consider the consensus alternating direction method of multipliers (ADMM) \cite{boyd2011distributed,goldstein2016unwrapping,huang2016consensus} as well as its asynchronous variant (async-ADMM), which aims at reducing latencies and thereby reduce the time-to-solution \cite{zhang2014asynchronous}. Consensus ADMM has previously been applied to high-dimensional inverse problems in data sciences \cite{miao2013hypergraph,ma2018fast}, statistical learning \cite{boyd2011distributed,xu2017adaptive,goldstein2016unwrapping,AAAIW1715174}, and imaging \cite{goldstein2014fast,juesas2015consensus,heredia2017norm}. The algorithm tackles large-scale problems by partitioning the data into, say, $N$ smaller batches that can be solved in parallel, and in some cases explicitly.
This often leads to an improved ratio of local computation and communication. 
 More specifically, each iteration of the algorithm breaks down into $N$ subproblems using parts of the data, an averaging step that is performed once their corresponding processors have solved all $N$ subproblems, and an explicit update of the dual variable. The main change in the async-ADMM variant is that the averaging step is performed once $N_a<N$ subproblems have been solved, reducing the overall latency.

As we demonstrate in our numerical experiments, a straightforward implementation of consensus ADMM converges slowly in particular when the information contained in the split data sets is complementary and the number of batches, $N$, is large. One problem in these cases is that the averaging step in consensus ADMM gives equal weighting to all the solutions corresponding to each batch, leading to an uninformed averaged reconstruction. 
In large-scale problems such as PDE parameter estimation, this renders consensus ADMM prohibitive since often only a few iterations are affordable.

To increase the performance of consensus ADMM, particularly in early iterations, we introduce a novel weighting scheme that improves the convergence of consensus ADMM. The weights are obtained in a systematic and efficient way using the framework of uncertainty quantification (UQ) proposed in \cite{flath2011fast}. We demonstrate the effect of the weights on a collection of linear inverse problems. We also outline the potential of our method by comparing it to the Gauss-Newton method \cite{haber2014computational} and the Nonlinear CG method \cite{hager2005new} on a single-physics PDE parameter estimation problem involving a travel time tomography survey, and a multiphysics parameter estimation problem involving Direct Current Resistivity (DCR) and travel time tomography \cite{treister2016fast} surveys. 

The remainder of the paper is organized as follows. In Sec.~\ref{sec:MathBackground}, we review the mathematical framework for MAP estimation and UQ as well as numerical optimization algorithms for their computations. In Sec.~\ref{sec:wADMM}, we present the weighted consensus ADMM method and its asynchronous variant. In Sec.~\ref{sec:num},  we outline the potential of our method for with a series of numerical experiments, and finally, we summarize the paper in Sec.~\ref{sec:conclusion}.

\section{Mathematical Background and Numerical Implementation} \label{sec:MathBackground}
In this section, we briefly review the computation of the MAP estimate and principles of uncertainty quantification in the context of large-scale Bayesian inverse problems (see, e.g., \cite{calvetti2006large,stuart2010inverse} for a detailed overview). We limit the discussion to the finite-dimensional case since we follow the discretize-optimize approach, however, an overview of the infinite-dimensional case can be found in \cite{stuart2010inverse}. We also review optimization techniques for computing the MAP estimate and their associated challenges. 

\subsection{MAP Estimation and UQ}
We consider additive noise-corrupted observations
\begin{align}
  Y = \mathcal{F}(X) + E,
\end{align}
where $\mathcal{F} \colon \bbR^n \mapsto \bbR^{m}$ is the parameter-to-observable map, and $Y, X,$ and $E$ are random variables corresponding to the observations, the model parameter , and the measurement noise, respectively. In the following, we denote observations by $\bfy \in \bbR^{m}$,  $\bfx \in \bbR^n$, and $\bfepsilon \in \bbR^m$.

We employ the prior probability distribution function (PDF), $\pi_{\rm prior} \colon \bbR^n \mapsto \bbR$, which describes prior information we may have about $X$, and the likelihood function $\pi(\bfy | \bfx)$ which describes the relationship between the measurements $\bfy$ and the unknown model parameters $\bfx$. We use Bayes' Theorem to obtain the posterior PDF
 \begin{align}
    \pi_{\rm post}(\bfx) \propto \pi_{\rm prior}(\bfx) \pi(\bfy | \bfx),
    \label{eq:Bayes}
 \end{align}
 and compute the MAP point by maximizing the posterior distribution, that is,
 \begin{align}
    \bfx_{\rm MAP} = \argmax_{\bfx} \pi_{\rm post}(\bfx).
 \end{align}

For simplicity, we assume that $X$ and $E$ are statistically independent
and limit the discussion to the case where the prior PDF is Gaussian and $E$ is independently and identically distributed, i.e., $E \sim \mathcal{N}({\bf0}, \Gamma_{\rm noise}),$ where $\Gamma_{\rm noise} \in \bbR^{m \times m}$ is the diagonal noise-covariance matrix. In this case, the likelihood and prior PDFs are given by 
\begin{align}
\label{eq:GaussianPDFs}
  \pi(\bfy|\bfx) \propto  \exp (-\Phi(\bfx))
  \quad \text{ and } \quad 
  \pi_{\rm prior}(\bfx) \propto \exp (-\mathcal{R}(\bfx)),
\end{align}
respectively, where due to the assumptions above
\begin{align}\label{eq:misfitReg}
  \Phi(\bfx) = \hf \| \mathcal{F} (\bfx) - \bfy \|_{\boldsymbol{\Gamma}_{\rm noise}^{-1}}^2 \; \text{ and } \; \; \; 
  \mathcal{R}(\bfx) = \hf \|\bfx - \bfx_{\rm ref} \|_{\boldsymbol{\Gamma}_{\rm{prior}}^{-1}}^2.
\end{align}
Here, $\Phi, \mathcal{R}\colon \bbR^n \mapsto \bbR$ are the misfit and regularizer, respectively, $\bfx_{\rm ref}$ is the mean of the model parameter prior PDF, and $\Gamma_{\rm prior} \in \bbR^{n \times n}$ is the covariance matrix of the prior PDF. 
Using \eqref{eq:Bayes} and (\ref{eq:GaussianPDFs}), we can restate the posterior distribution in closed form as
\begin{align}
  \pi_{\rm post} \propto \exp \Big(-\Phi(\bfx) - \mathcal{R}(\bfx)\Big).
\end{align}
Finding the MAP estimate is then equivalent to solving the regularized weighted least squares problem
\begin{align}
  \label{eq:MAP}
  \bfx_{\rm MAP} = \argmin_{\bfx} \;\; \Big(\Phi(\bfx) + \mathcal{R}(\bfx)\Big).
\end{align}

When $\mathcal{F}$ is a linear operator, that is, $\mathcal{F} = \bfA \in \bbR^{m \times n}$, the posterior PDF is also Gaussian and we can write its covariance matrix $\Gamma_{\rm post} \in \bbR^{n \times n}$ in closed form as
\begin{align}
\label{eq:postCovIntractable}
  \boldsymbol{\Gamma}_{\rm post} = (\bfA^\top \boldsymbol{\Gamma}_{\rm noise}^{-1} \bfA + \boldsymbol{\Gamma}_{\rm prior}^{-1})^{-1},
\end{align}
which can be used for quantifying uncertainties of the model parameter $\bfx$. In the context of large-scale PDE parameter estimation, however, the matrix $\bfA$, let alone its inverse, is seldom constructed. The computation of $\bfGamma_{\rm{post}}$ is therefore intractable, which is why we follow~\cite{flath2011fast} and use an iterative method to obtain an approximation in Sec.~\ref{subsec:weights}.

\subsection{Numerical Optimization}
\label{subsec:parallelGaussNewton}
In large-scale PDE parameter estimation, we are concerned with the case where massive amounts of data are available, leading to numerous right-hand sides and potentially multiple PDEs \cite{treister2016fast, haber2013model}. Here, we split the misfit function in~\eqref{eq:misfitReg} into $N$ terms, i.e.,
\begin{equation}
  \Phi(\bfx) = \sum_{j=1}^N \Phi_j(\bfx), \quad\text{where}\quad \Phi_j(\bfx)   = \hf \| \mathcal{F}_j (\bfx) - \bfy_j \|_{\boldsymbol{\Gamma}_{j, \rm noise}^{-1}}^2,
\end{equation}
where $\mathcal{F}_j \colon \bbR^{n} \mapsto \bbR^{m_j}$ and $ \bfy_j \in \bbR^{m_j}$ correspond to the $j^{\text{th}}$ forward operator and right-hand side, respectively, and $\Gamma_{j, \rm noise} \in \bbR^{m_j \times m_j}$ is the noise covariance matrix corresponding to the $j^{\text{th}}$ misfit term.  We can rephrase (\ref{eq:MAP}) as
\begin{equation}
  \begin{split}
    \label{eq:MAPminProb}
    \bfx_{\rm MAP} = \argmin_\bfx \; \sum_{j=1}^N \Phi_j(\bfx) + \mathcal{R}(\bfx).
    \end{split}
\end{equation}
There many ways for exploiting the structure of~\eqref{eq:MAPminProb} including stochastic optimization methods, e.g., stochastic approximation~\cite{RobbinsMonro1951}, stochastic average approximation~\cite{KleywegtEtAl2006}, and the method of simultaneous sources~\cite{HaberChungHerrmann2012}.
Here, we are interested in deterministic optimization methods (potentially applied to a stochastic average approximation or the reduced problem in~\cite{HaberChungHerrmann2012}). Common choices include steepest descent (SD), quasi-Newton methods such as l-BFGS \cite{wright1999numerical}, Nonlinear Conjugate Gradient \break(NLCG) methods \cite{hager2006survey,hager2005new}, Gauss-Newton methods \cite{haber2014computational,wright1999numerical}. We limit the discussion in this section to the Gauss-Newton-PCG and NLCG methods for tackling large-scale PDE parameter estimation problem.

When applying the Gauss-Newton-PCG algorithm to~\eqref{eq:MAPminProb} the \emph{linearized} Hessian of the objective function is used as the coefficient matrix in the Newton system. The direction to update the model, $\partial \bfx \in \bbR^n$, is computed approximately by applying a preconditioned CG method (see, e.g.,~\cite{Saad2003}) to the linear system
\begin{align}
  \left( 
    \sum_{j=1}^{N} \bfJ_j^\top \bfGamma_{j,\rm{noise}}^{-1} \bfJ_j
    +
    \nabla_\bfx^2 \mathcal{R}(\bfx)
  \right)
  \partial \bfx = - \sum_{j=1}^{N} \nabla_{\bfx} \Phi_j(\bfx) - \nabla_\bfx \mathcal{R}(\bfx)
  \label{eq:GaussNewtonSystem}
\end{align}
(see Alg.~\ref{alg:GaussNewton}, step 2), where $\bfJ_j \in \bbR^{m_j \times n}$ is the Jacobian matrix of the $j^{\rm th}$ misfit function $\Phi_j$.
\begin{algorithm}[t]
  \begin{itemize}
    \item initialize $\bfx^{(0)}$
    \item for $k=1,2,\ldots$
    \begin{enumerate}
      \item compute $\Phi_1(\bfx^{(0)}),\ldots,\Phi_N(\bfx^{(0)})$ and $\nabla_{\bfx} \Phi_1(\bfx^{(0)}),\ldots, \nabla_{\bfx} \Phi_N(\bfx^{(0)})$
      \item solve system~\eqref{eq:GaussNewtonSystem} to obtain $\delta \bfx$ using PCG
      \item set $\bfx^{(k+1)} = \bfx^{(k)} + \gamma \partial \bfx$, where $\gamma$ is set by a linesearch
      \item check convergence criteria
    \end{enumerate}
  \end{itemize}
 \caption{Gauss-Newton}
 \label{alg:GaussNewton}
\end{algorithm}
\begin{algorithm}[t]
  \begin{itemize}
    \item initialize $\bfx^{(0)}$
    \item set $\bfp^{(0)} = - \nabla_\bfx f(\bfx^{(0)}),$ where $f(\bfx) = \sum_{j=1}^N \Phi_j(\bfx) + \mathcal{R}(\bfx)$
    \item for $k=1,2,\ldots$
    \begin{enumerate}
      \item update $\bfx^{(k+1)} \leftarrow \bfx^{(k)} + \gamma \bfp^{(k)}$, where $\gamma$ is set by a linesearch
      \item compute $\nabla_\bfx f(\bfx^{(k+1)})$
      \item set $\bfd^{(k)} = \nabla_\bfx f(\bfx^{(k+1)}) - \nabla_\bfx f(\bfx^{(k)})$
      \item compute $\beta^{(k)} = \dfrac{1}{(\bfp^{(k)})^\top \bfd^{(k)}} \left(\bfd^{(k)} - 2\bfp^{(k)} \dfrac{\| \bfd^{(k)} \|^2}{(\bfp^{(k)})^\top \bfd^{(k)}}\right) \nabla_\bfx f(\bfx^{(k+1)})$
      \item update $\bfp^{(k+1)} = -\nabla_\bfx f(\bfx^{(k+1)}) + \beta^{(k)}\bfp^{(k)}$
      \item check convergence criteria
    \end{enumerate}
  \end{itemize}
 \caption{NLCG}
 \label{alg:NLCG}
\end{algorithm}
Although the individual terms in~\eqref{eq:GaussNewtonSystem} can be computed in parallel, an efficient implementation is non-trivial. To limit the communication overhead, one can use  the static scheduling approach described in~\cite{ruthotto2017jinv}.
Here, the model and a number of meshes, sources, receivers, and forward problems are assigned to all the workers in the offline phase. Then, to evaluate the misfit, compute one full gradient, and perform matrix-vector products with the Hessian of the objective function \eqref{eq:MAPminProb}, each worker computes its corresponding batch of gradients and Hessians and communicates it to the main process. Note that each matrix-vector product with the Hessian, which is done in each PCG iteration for solving~\eqref{eq:GaussNewtonSystem}, requires sending and receiving a vector to and from each worker. For large-scale problems this can result in a nontrivial amount of communication, especially when many PCG iterations are needed. Moreover, if the data is  divided unevenly among the workers, the algorithm may lead to large latencies in each PCG iteration. This motivates us to consider more scalable distributed algorithms, especially when the size and dimension of the problem is very large.

The NLCG algorithm~\cite{hager2005new} requires substantially less communication per outer iteration. NLCG performs explicit steps using gradients to update the model (see Alg.~\ref{alg:NLCG}), and therefore avoids the communication that comes with solving the Gauss-Newton system using an iterative method. In our experience, however, the method requires many more iterations than the Gauss-Newton-PCG method in order to achieve the same level of accuracy (see Sec.~\ref{subsec:singlePhysics} and \ref{subsec:jointInversion}). Since each gradient and objective function evaluation requires at least $N$ PDE solves, the large number of outer iterations renders the NLCG method less attractive for solving these large-scale PDE parameter estimation problems. 
\section{Uncertainty-Weighted Consensus ADMM} \label{sec:wADMM}
In this section, we introduce our uncertainty-weighted ADMM method. First, we present the general formulation of the weighted ADMM, which involves rephrasing (\ref{eq:MAPminProb}) as a global variable consensus problem \cite{boyd2011distributed}, and review the asynchronous implementation presented in~\cite{zhang2014asynchronous}. Following \cite{flath2011fast}, we use a systematic scheme for selecting the weights, which is based on the uncertainties of the local models. Finally, we use a numerical example to illustrate the intuition behind the weights.
\subsection{Weighted Consensus ADMM} 
Motivated by the discussion in the previous section, we reformulate the optimization problem (\ref{eq:MAPminProb}) as an equivalent weighted global variable consensus problem
\begin{equation}
\label{eq:globalVarConsMAP}
\begin{split}
  \bfx_{\rm MAP} = \: &\argmin_{\bfx_1,\ldots,\bfx_N,\bfz} \; \sum_{j=1}^N \left(\Phi_j(\bfx_j) + \mathcal{R}(\bfx_j)\right), \\ 
  &\text{ s.t. } \quad \bfW_j(\bfx_j-\bfz) = {\bf0}, \quad j=1,\ldots,N,
\end{split}
\end{equation}
where in contrast to (\ref{eq:MAPminProb}), the objective function is now separable, however, the coupling is enforced by the constraints.
Here, $\bfx_j \in \bbR^n$ are the local variables that are brought into consensus via the global variable $\bfz \in \bbR^n$, and $\bfW_j \in \bbR^{n \times n}$ are diagonal weight matrices. In the standard global consensus formulation \cite{boyd2011distributed}, the identity matrix is assigned as the weight matrices. This reformulation allows each of the objective terms in (\ref{eq:globalVarConsMAP}) to be handled by its corresponding worker via the consensus ADMM algorithm.

Consensus ADMM aims at solving problem (\ref{eq:globalVarConsMAP}) by finding a saddle point of its Lagrangian via the following iterations:
\begin{align}
        \bfx_j^{(k+1)} &= \argmin_{\bfx_j} \; \left( \Phi_j(\bfx_j) + \mathcal{R}(\bfx_j) + (\bfu_j^{(k)})^\top \bfW_j\bfx_j + \frac{\rho}{2}\| \bfW_j(\bfx_j - \bfz^{(k)}) \|_2^2 \right), \label{alg:xstepExplicit} \\ j=1,&\ldots,N, \nonumber
    \\
    \bfz^{(k+1)} &= \bigg(\sum_{j=1}^N \bfW_j^\top \bfW_j\bigg)^{-1} \sum_{j=1}^N \left( \bfW_j^\top \bfW_j\bfx_j^{(k+1)} + (1/\rho) \bfW_j\bfu_j^{(k)} \right), \label{alg:zstepExplicit}
    \\
    \bfu_j^{(k+1)} &= \bfu_j^{(k)} + \rho \bfW_j(\bfx_j^{(k+1)} - \bfz^{(k+1)}), \quad j=1,\ldots,N,\label{alg:ustepExplicit}
\end{align}
where $k$ denotes the current iteration, $\bfu_j$ are the dual variables, and $\rho>0$ is the penalty parameter associated with the augmented Lagrangian term. 

The minimization steps in~\eqref{alg:xstepExplicit} entail PDE solves and are the most computationally challenging part of the algorithm, however, they correspond to the local subproblems that are solved independently by each processor.
A further advantage is that the local subproblem can be solved using any optimization algorithm, which provides an easy way to tailor the method to different subproblems, e.g., subproblems containing different PDEs for which highly-optimized algorithms already exist. 
Consequently, ADMM sits at a higher-level of abstraction from classical optimization algorithms such as those mentioned in Sec.~\ref{subsec:parallelGaussNewton}.
The global variable $\bfz$ attempts to bring the local variables $\bfx_j$ into consensus by averaging them in (\ref{alg:zstepExplicit}), and finally, the dual variables are updated via a gradient ascent step in (\ref{alg:ustepExplicit}). 

We use the stopping criteria in~\cite{boyd2011distributed} to define the primal and dual residuals as 
\begin{align}
\bfr^{(k+1)} &= \left(\bfW_1 (\bfx_1^{(k+1)} - \bfz^{(k+1)}),\ldots, \bfW_N(\bfx_N^{(k+1)}-\bfz^{(k+1)})\right), \; \label{eq:primalRes}\text{ and} \\
\bfs^{(k+1)} &= -\rho\left(\bfW_1(\bfz^{(k+1)} - \bfz^{(k)}),\ldots,\bfW_N(\bfz^{(k+1)} - \bfz^{(k)})\right), \label{eq:dualRes}
\end{align}
respectively, and stop whenever 
\begin{align}
\label{eq:stopCrit}
  \| \bfr^{(k)} \|_2 \leq \epsilon^{\rm pri} \; \text{ and } \; \| \bfs^{(k)} \|_2 \leq \epsilon^{\rm dual}
\end{align}
for some chosen primal and dual tolerances $\epsilon^{\rm pri}$ and $\epsilon^{\rm dual}$. It is also common to adaptively choose the penalty parameter $\rho$. We use the scheme in~\cite{boyd2011distributed}, i.e.,
\begin{equation} 
\label{eq:adaptiveRho}
\rho^{(k+1)} =   \left\{
\begin{array}{ll}
      \tau^{\rm incr} \rho^{(k)} & \text{if} \quad \Vert \bfr^{(k)} \Vert_2 > \mu \Vert \bfs^{(k)} \Vert_2 \\
      \rho^{(k)}/ \tau^{\rm decr} & \text{if} \quad \Vert \bfs^{(k)} \Vert_2 > \mu \Vert \bfr^{(k)} \Vert_2 \\
      \rho^{(k)} & \text{otherwise},
\end{array} 
\right. 
\end{equation}
where $\mu>1, \tau^{\rm incr}>1, $ and $\tau^{\rm decr}>1$ are parameters commonly chosen to be $10, 2,$ and $2$, respectively \cite{boyd2011distributed}. This updating scheme aims at balancing the primal and dual residual norms within a factor of $\mu$ of one another as they both converge to zero.

Parallelization of consensus ADMM is much more straightforward than that of the Gauss-Newton-PCG described in Sec.~\ref{subsec:parallelGaussNewton}. The amount of communication per outer iteration is reduced as we only communicate one set of models, $\bfx_1,\ldots,\bfx_N$. In the synchronous parallel implementation, the master processor must wait for all the workers to finish solving their corresponding subproblems in \eqref{alg:xstepExplicit} before performing the averaging step \eqref{alg:zstepExplicit} per iteration, which may lead to high latencies when some of the workers are much slower than others. The async-ADMM method in~\cite{zhang2014asynchronous} aims at reducing latencies in star network topologies. Here, the global averaging step \eqref{alg:zstepExplicit} is performed when $N_a < N$ workers report their results. A \textit{bounded delay} condition is also enforced, where every worker has to report at least once every $k_a$ iterations to ensure sufficient "freshness" of all updates. We note that here we have better control of the overall amount of communication and latency since we can administer how many forward problems to assign to any given worker, and how accurately to solve each subproblem.

Convergence results have been established for the synchronous ADMM algorithm in the case where the local subproblems are convex. In this case, the algorithm converges regardless of the initial choice $\rho^{(0)}$ \cite{hong2017linear,eckstein1992douglas}.
Even when solving \eqref{alg:xstepExplicit} inexactly, ADMM convergence can be shown~\cite[Sec. 4]{hong2017linear}.
For the asynchronous case, convergence is ensured via the bounded delay condition. However, additional factors such as network bandwidth and processor configuration are taken into account \cite{zhang2014asynchronous}. For non-convex subproblems, it has been shown that ADMM converges to a local minimum under some modest assumptions, most importantly requiring $\rho$ to be sufficiently large, \cite{Macdonald2018,wright1999numerical,hong2016convergence}. These assumptions ensure that the Hessian of the Lagrangian of (\ref{eq:globalVarConsMAP}) remains positive definite throughout the ADMM iterations.
\begin{algorithm}[t]
  \begin{itemize}
    \item initialize $\bfx_j^{(0)}, \bfz^{(0)},$ and $\bfu_j^{(0)}$ for $  j=1,\ldots,N$
    \item while~\eqref{eq:stopCrit} not satisfied
    \begin{enumerate}
      \item solve \eqref{alg:xstepExplicit} locally
      \item update $\bfz$ using the averaging step \eqref{alg:zstepExplicit}
      \item update dual variables \eqref{alg:ustepExplicit}
    \end{enumerate}
  \end{itemize}
 \caption{Consensus ADMM}
 \label{alg:wADMM}
\end{algorithm}
\begin{algorithm}[t]
  \begin{itemize}
    \item initialize $\bfx_j^{(0)}, \bfz^{(0)},$ and $\bfu_j^{(0)}$ for $ j=1,\ldots,N$
    \item initialize $N_a$ and $k_a$
    \item while~\eqref{eq:stopCrit} not satisfied
    \begin{enumerate}
      \item solve \eqref{alg:xstepExplicit} locally
      \item perform averaging step \eqref{alg:zstepExplicit} when $N_a$ workers report their solutions
      \item update the corresponding $N_a$ dual variables \eqref{alg:ustepExplicit}
    \end{enumerate}
  \end{itemize}
 \caption{Consensus async-ADMM}
 \label{alg:async-wADMM}
\end{algorithm}

\subsection{Selection of Weights}\label{subsec:weights}
We choose the weights to be the inverse of the diagonals of the posterior covariance $\Gamma_{j, \rm post} \in \bbR^{n \times n}$ corresponding to the $j^{\text{th}}$ objective term in (\ref{eq:globalVarConsMAP}). This is one way to assign higher weights to elements of $\bfx_j$ for which the $j^{\text{th}}$ subproblem contains more information. It also reduces the impact of elements for which the data of the subproblem is uninformative. Clearly, there are other options to transform uncertainties into weights. Since we are mostly interested in encoding large differences in the uncertainties between subproblems, we do not compute the uncertainties with high accuracy. 

As seen in (\ref{eq:postCovIntractable}), construction of the posterior covariance may not be tractable, especially for large-scale PDE parameter estimation problems and when the forward model is nonlinear. As a result, we follow the works of \cite{flath2011fast} for approximating the posterior covariance of each objective term in a tractable way. This is done via a low-rank approximation of the approximate Hessian of the misfit $\Phi_j$ in the following manner:
\begin{enumerate}
  \item We linearize the residual in $\Phi_j$ and obtain the Gauss-Newton approximation
  \begin{align}
  \bfH_{j, \rm mis} \approx \bfJ_j^\top (\boldsymbol{\Gamma}_{j, \rm noise}^{-1}) \bfJ_j,
  \end{align} 
  where $\bfJ_j \in \bbR^{m_j \times n}$ is the Jacobian matrix of $\mathcal{F}_j$ evaluated at some reference model parameter, e.g., $\bfx_{\rm ref}$. We note that explicit construction of $\bfH_{j, \rm mis}$ is not necessary as we only need the action of $\bfJ_j$ and $\bfJ_j^\top$ on a vector.
  \item Denoting the prior-conditioned approximate Hessian by $$\tilde{\bfH}_{j, \rm mis} = \boldsymbol{\Gamma}^{1/2}_{\rm prior} \bfH_{j, \rm mis} \boldsymbol{\Gamma}^{1/2}_{\rm prior},$$ we rewrite the $j^{\rm th}$ posterior covariance in (\ref{eq:postCovIntractable}) as
    \begin{equation}
    \label{eq:tempGammaPost}
      \boldsymbol{\Gamma}_{j, \rm post} = \boldsymbol{\Gamma}_{\rm prior}^{1/2} \left(\tilde{\bfH}_{j, \rm mis} + \bfI\right)^{-1} \boldsymbol{\Gamma}_{\rm prior}^{1/2}.
    \end{equation}
  \item We then construct a low-rank approximation of the prior-conditioned Hessian using, e.g., randomized SVD \cite{saibaba2016randomized} or Lanczos bidiagonalization \cite{golub2012matrix} to obtain
    \begin{equation}
      \tilde{\bfH}_{j, \rm mis} = \bfV \boldsymbol{\Lambda} \bfV^\top \approx \bfV_r \boldsymbol{\Lambda}_r \bfV_r^\top,
    \end{equation}
    where $\boldsymbol{\Lambda} = \text{diag}(\lambda_1,\ldots,\lambda_n) \in \bbR^{n \times n}$ and $\bfV = [ \bfv_1,\ldots, \bfv_n ] \in \bbR^{n \times n}$ denote the matrix of eigenvalues and eigenvectors of $\tilde{\bfH}_{j, \rm mis}$, respectively, and $\boldsymbol{\Lambda}_r = \text{diag}(\lambda_1,\ldots,\lambda_r) \in \bbR^{r \times r}$ and $\bfV_r =  [ \bfv_1,\ldots, \bfv_r ] \in \bbR^{n \times r}$ are their corresponding truncations retaining only the $r$ largest eigenvalues and eigenvectors.
  \item We plug this approximation into (\ref{eq:tempGammaPost}) and use the Sherman-Morrison-\\Woodbury formula \cite{sherman1950adjustment} to obtain an expression for the inverse term:
    \begin{equation}
      \left(\tilde{\bfH}_{j, \rm mis} + \bfI\right)^{-1} \approx \bfI - \bfV_r \bfD_r \bfV_r^\top + \mathcal{O}\left(\sum_{i=r+1}^n \frac{\lambda_i}{\lambda_i+1}\right),
    \end{equation}
    where $\bfD \in \bbR^{r \times r} = \text{diag}(\lambda_1/(\lambda_1 + 1),\ldots,\lambda_r/(\lambda_r + 1))$.
  \item Finally, we obtain a manageable approximation of the posterior covariance that does not involve any inverse terms:
    \begin{equation}
      \boldsymbol{\Gamma}_{j, \rm post} \approx \bfGamma_{\rm prior}^{1/2} (\bfI - \bfV_r \bfD_r \bfV_r^\top)\bfGamma_{\rm prior}^{1/2}.
    \end{equation}
\end{enumerate}

We choose the weights to be the inverse of the diagonals of $\bfGamma_{j, \rm post}$, 
\begin{equation}
  \bfW_j = \text{diag}(\bfGamma_{j, \rm post})^{-1}, \quad j=1,\ldots,N,
\end{equation}
so that we get higher weights in parts of the model where we are more certain and vice-versa. We may also update the weights in every iteration of our optimization scheme so that we instead employ local approximations of our posterior PDF \cite{bui2013computational}. In our experience, however, computing the weights once in the offline phase is enough to accelerate convergence in the early iterations. We note that this is only one way to estimate the diagonal entries of the posterior covariance matrix, and that a plethora of alternative options can be used, e.g., probing methods \cite{tang2012probing}, extrapolation methods \cite{fika2018estimating}, stochastic methods \cite{bekas2007estimator}, and domain decomposition methods \cite{lin2009fast,tang2011domain,li2008computing}. 
  \begin{figure}[t]
    \centering
    \begin{tabular}{cccc}
      $\bfx_1^{(1)}$ & $\bfx_2^{(1)}$ & & $\bfz^{(1)}$ unweighted
      \\
      \includegraphics[width=0.22\textwidth, height=1.35in]{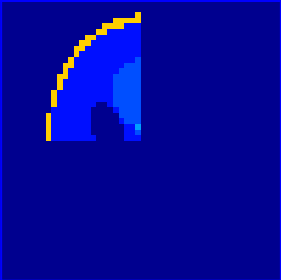}
      &
      \includegraphics[width=0.22\textwidth, height=1.35in]{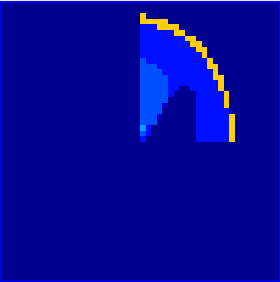}
      &
      &
      \includegraphics[width=0.22\textwidth, height=1.35in]{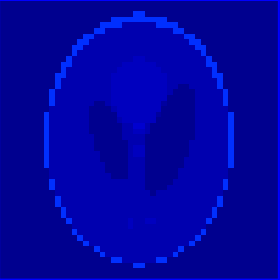}
      \\
      $\bfx_3^{(1)}$ & $\bfx_4^{(1)}$ & & $\bfz^{(1)}$ weighted
      \\
      \includegraphics[width=0.22\textwidth, height=1.35in]{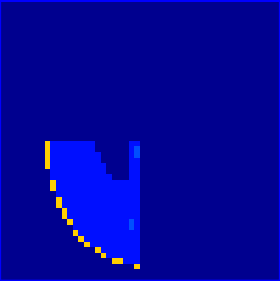}
      &
      \includegraphics[width=0.22\textwidth, height=1.35in]{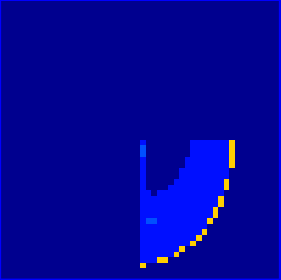}
      &
      &
      \includegraphics[width=0.22\textwidth, height=1.35in]{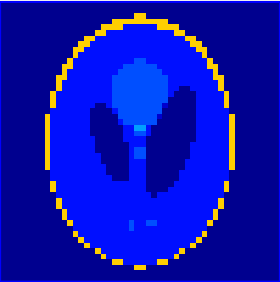}
    \end{tabular}
    \caption{\textit{Averaging step of the weighted and unweighted consensus ADMM for Ex. \ref{ex:identity}. In this case, $\bfW_1$ assigns higher weights to the pixels in the upper-left quadrant of $\bfx_1$, $\bfW_2$ assigns higher weights to the upper-right quadrant of $\bfx_2$, etc. As a result, the weights educate the averaging step, leading to a better reconstruction of the image.}}
    \label{fig:identityExampleWeighted}
    \vspace{-5mm}
\end{figure}
When the weights equal to one, the weighted ADMM method corresponds to the standard unweighted ADMM scheme, which is known to have slow convergence \cite{boyd2011distributed}. One reason is that the averaging step in (\ref{alg:zstepExplicit}) gives equal weighting to all elements of $\bfx_j$ for all $j=1,\ldots,N$, leading to poor reconstructions of $\bfz$, especially in the early iterations. To illustrate this, we perform the following example.
\begin{example} \label{ex:identity}
  Consider solving the trivial linear system $\bfI \bfx = \bfy$ with the weighted and unweighted consensus ADMM with $N=4$ splittings, where $\bfI \in \bbR^{n \times n}$ is the identity matrix, and $\bfx, \bfy \in \bbR^n$ are the model and the observed data, respectively. We formulate the least squares problem as 
  \begin{align}
    &\argmin_{\bfx_j, \bfz} \;\; \sum_{j=1}^4 \left( \hf \| \bfI_j \bfx_j - \bfy_j \|_2^2 + \frac{\alpha}{2} \| \bfL \bfx_j \|_2^2 \right)
    \\
    &\;\;\text{\normalfont  s.t. } \quad \bfW_j(\bfx_j - \bfz) = \boldsymbol{0}, \;\;\; j=1,\ldots,4,
  \end{align}
  where $\bfI_j \in \bbR^{(n/4) \times n}$ and $\bfy_j \in \bbR^{n/4}$ are subsets of the data obtained by partitioning the rows of $\bfI$ and $\bfy$ corresponding to the pixels in the top left, top right, bottom left, and bottom right quadrant of the domain as seen in Fig~\ref{fig:identityExampleWeighted}. We show the averaged reconstruction of both methods during the first iteration in Fig~\ref{fig:identityExampleWeighted}.
  
The weights associated with the subproblems, in this case, look very different from one another (see Fig.~\ref{fig:identityExampleWeighted}). Thus, introducing the weights considerably improves the effectiveness of the averaging step in ADMM and leads to faster convergence, especially in the early iterations. Our intuition is that when the weights look similar, the weighted ADMM will perform a similar averaging to that of the unweighted ADMM. This leads to a comparable performance of both methods. We illustrate this in the tomography problem in Sec.~\ref{subsec:LeastSquares}. 
\end{example}

\section{Numerical Experiments}\label{sec:num}
In this section, we outline the potential of the weighted scheme for consensus ADMM as well as its asynchronous variant on a series of linear and nonlinear inverse problems. We first experiment on a deblurring and a tomography problem from \texttt{Regtools}, a MATLAB package containing discrete ill-posed inverse problems~\cite{hansen1994regularization}, as well as from a collection of linear least-squares problems from the UF Sparse Matrix Collection \cite{davis2011university}. We then test our method on larger 3D PDE parameter estimation problems: a single-physics parameter estimation problem involving a travel-time tomography survey, and a multiphysics parameter estimation problem involving DCR and travel-time tomography.
\subsection{Least-Squares} \label{subsec:LeastSquares}
\newcommand{\rottext}[1]{\rotatebox{90}{\hbox to 30mm{\hss #1\hss}}}
\newcommand{\rottextt}[1]{\rotatebox{90}{\hbox to 25mm{\hss #1\hss}}}
\begin{figure}[t]
  \small
  \centering
  \begin{tabular}{cccc}
    & \textbf{ground truth} & \textbf{weighted ADMM} & \textbf{unweighted ADMM}
    \\
    \rottext{\small{\textbf{deblurring}}} \hspace{-4mm}
    &
    \includegraphics[width=0.28\textwidth,height=1.35in]{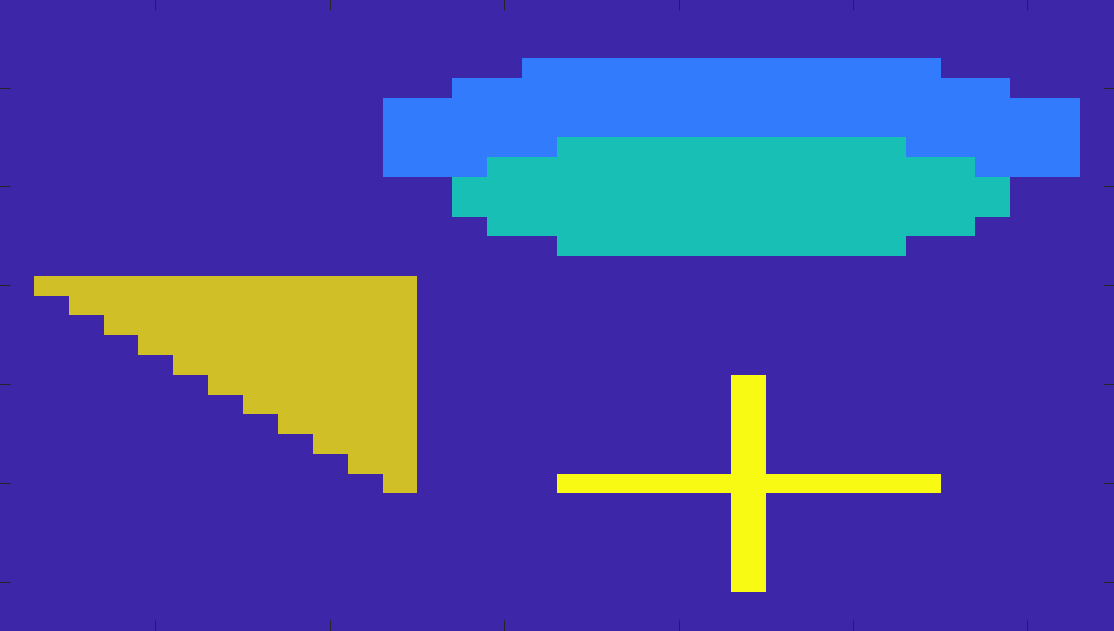}
    &
    \includegraphics[width=0.28\textwidth,height=1.35in]{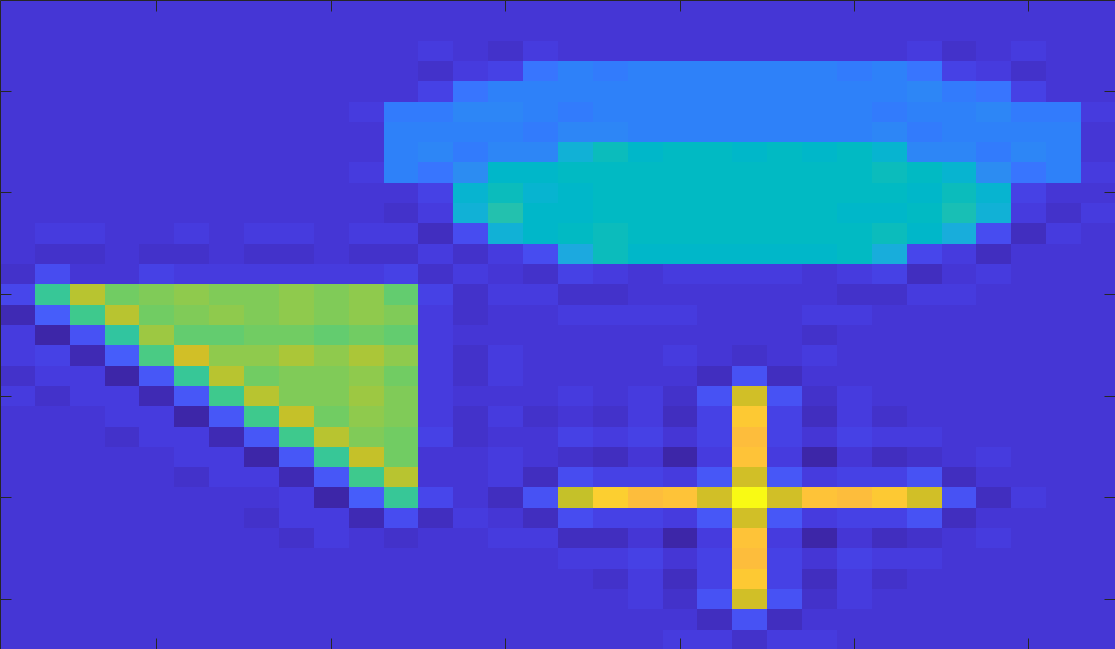}
    &
    \includegraphics[width=0.28\textwidth,height=1.35in]{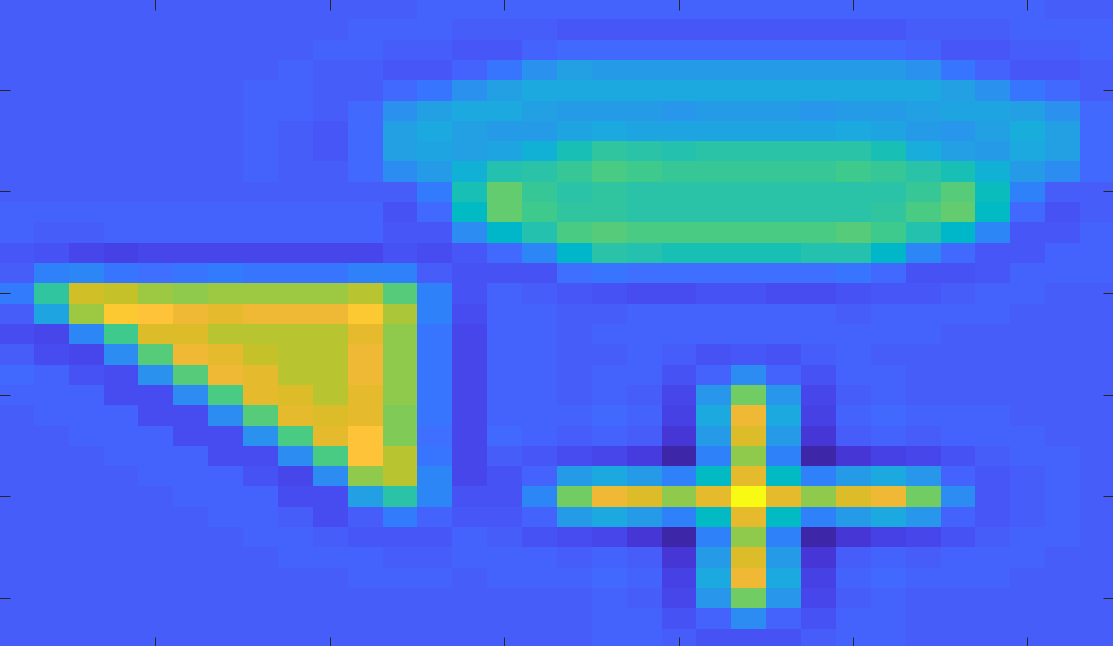}
    \\
    & & \textbf{relerr} $\approx $ 1.27e-01 & \textbf{relerr} $\approx$ 2.99e-01
  \end{tabular}
  \begin{tabular}{ccccc}
    & $\bfW_1$ & $\bfW_2$ & $\bfW_3$ & $\bfW_4$ \hspace{-5mm}
    \\
    \rottextt{\small{\textbf{Weights}}} \hspace{-5mm}
    &
    \includegraphics[width=0.205\textwidth,height=1.2in]{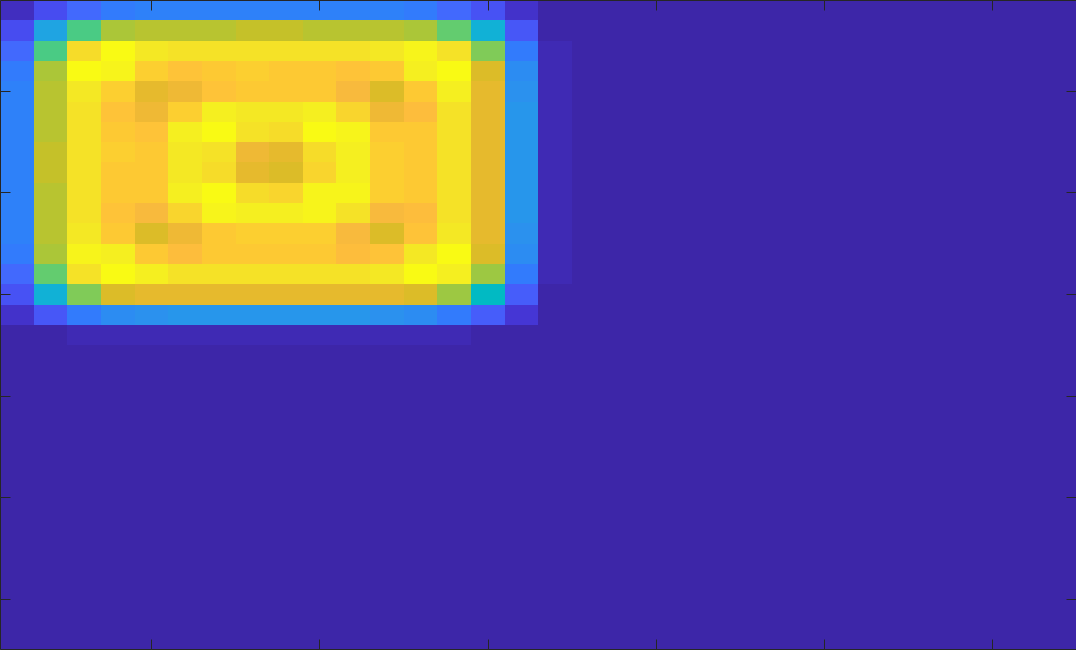}
    &
    \includegraphics[width=0.205\textwidth,height=1.2in]{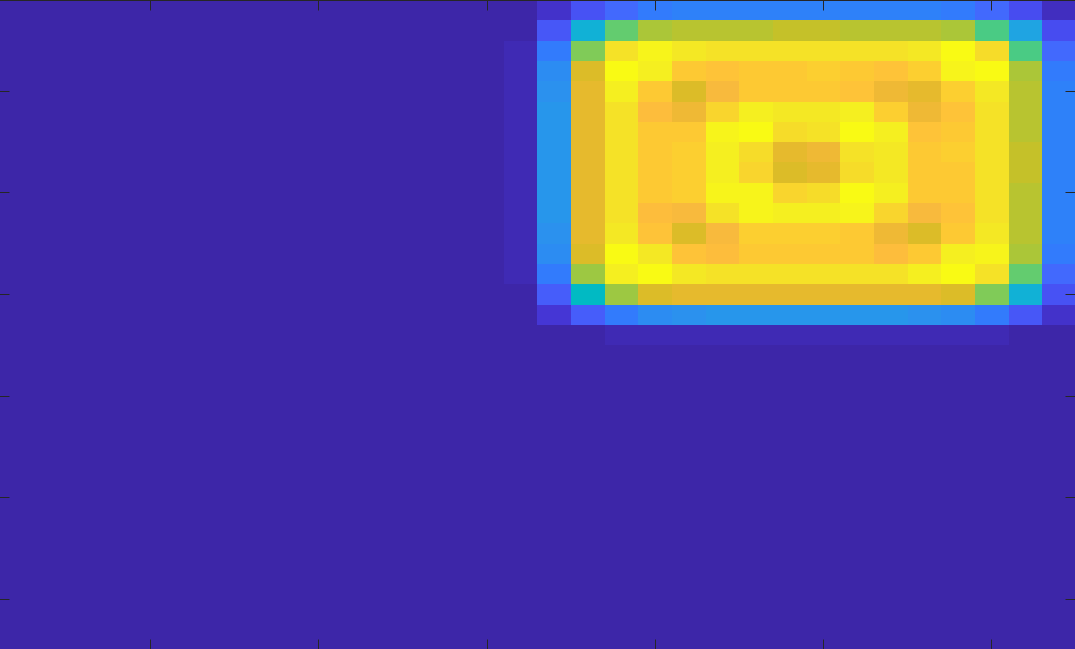}
    &
    \includegraphics[width=0.205\textwidth,height=1.2in]{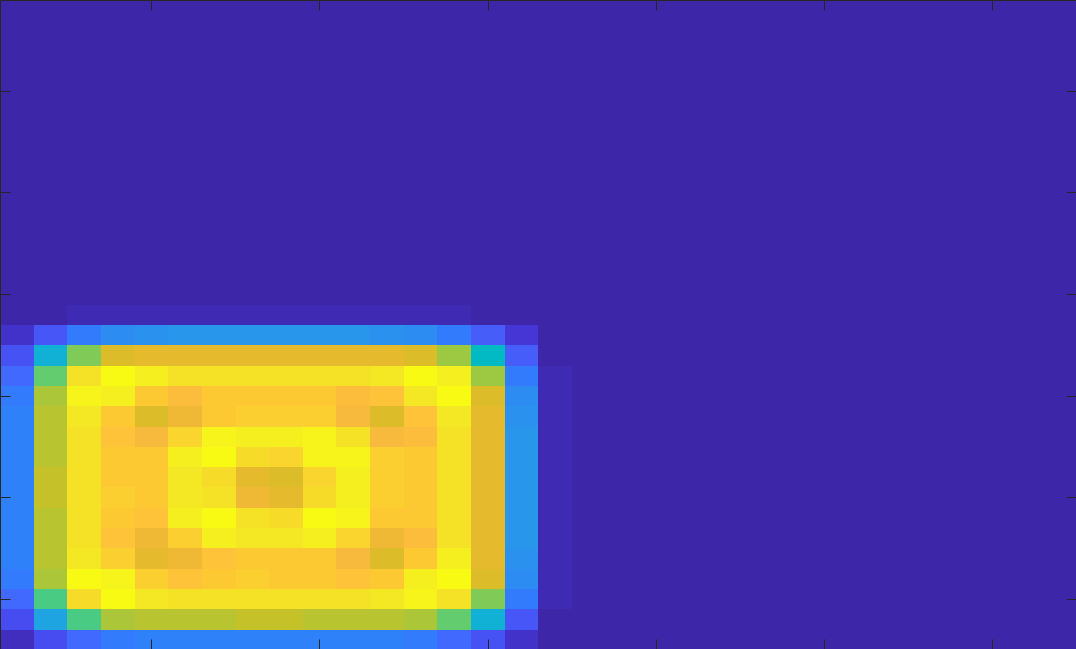}
    &
    \includegraphics[width=0.205\textwidth,height=1.2in]{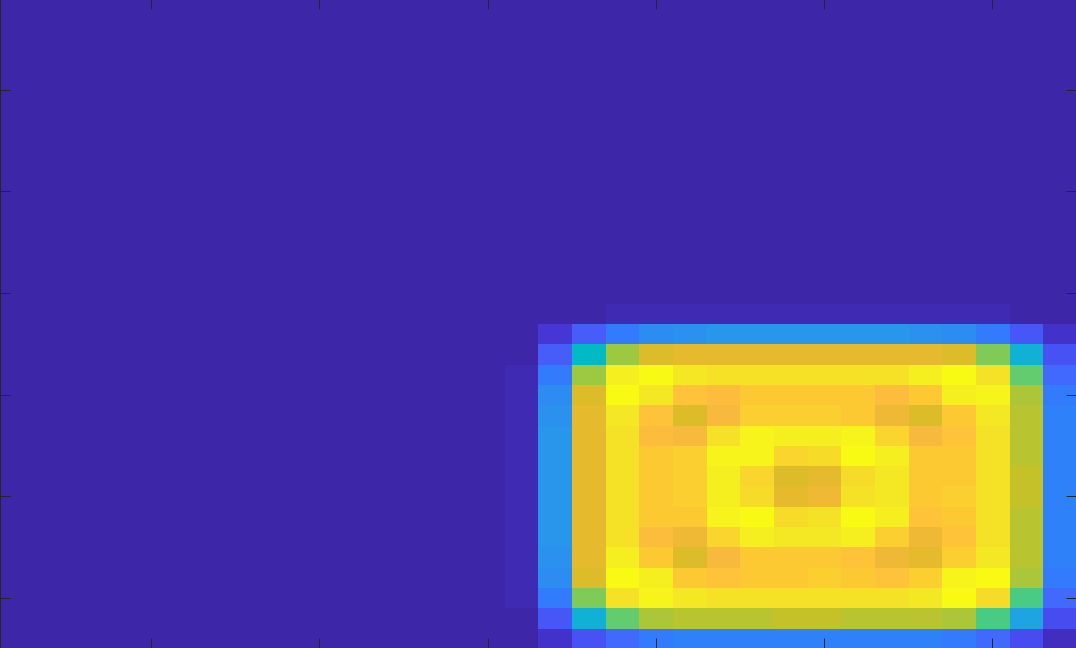}   
  \end{tabular}
  \caption{\textit{Reconstructions (first row) after 10 iterations and weights (second row) for the deblurring problem from \emph{\texttt{Regtools}~\cite{hansen1994regularization}}. 
  }}
  \label{fig:blur}
\end{figure}
\begin{figure}[!t]
  \small
  \centering
  \begin{tabular}{cccc}
    & \textbf{ground truth} & \textbf{weighted ADMM} & \textbf{unweighted ADMM}
    \\
    \rottext{\small{\textbf{tomography}}} \hspace{-4mm}
    &
    \includegraphics[width=0.28\textwidth,height=1.35in]{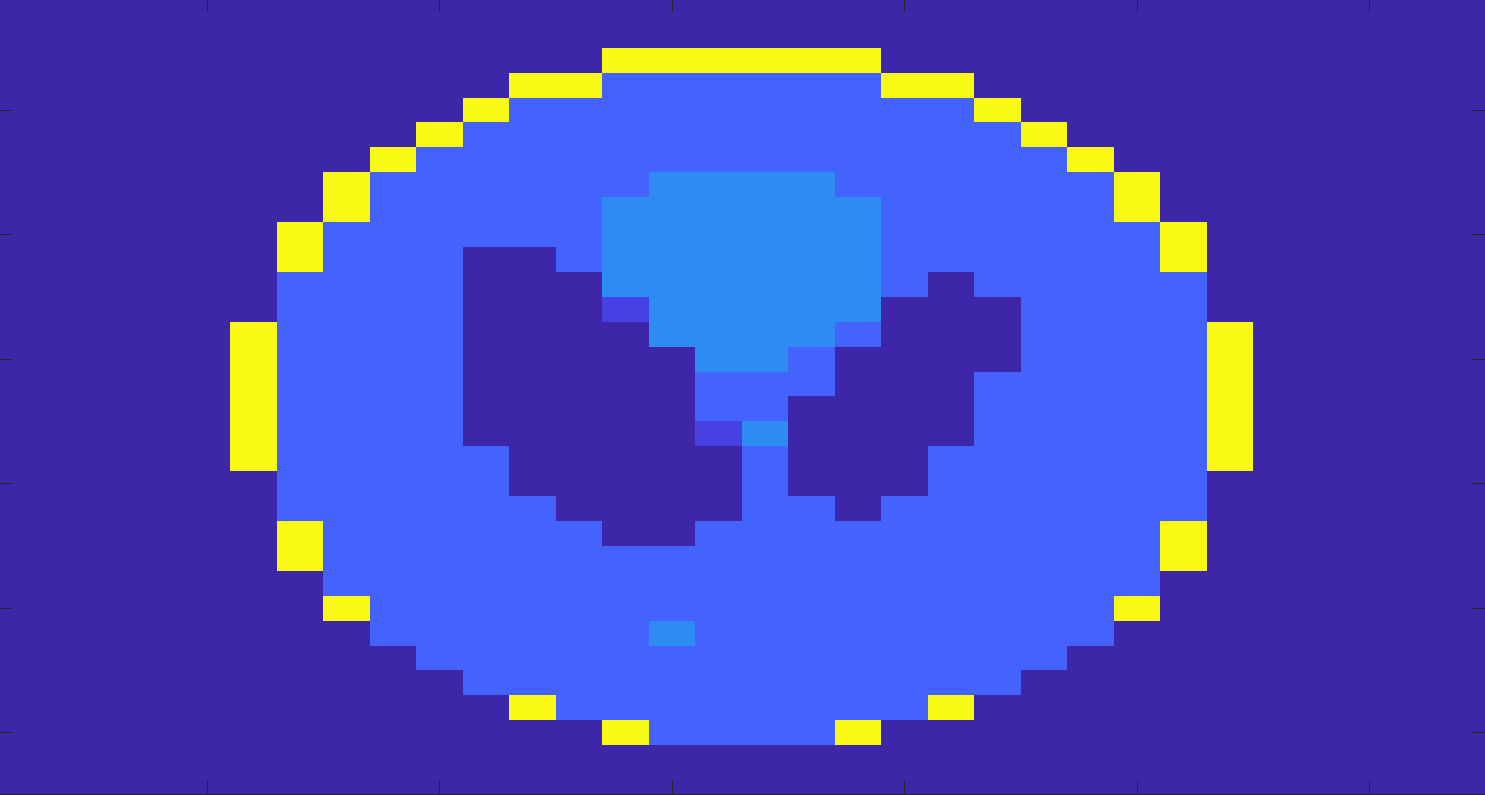}
    &
    \includegraphics[width=0.28\textwidth,height=1.35in]{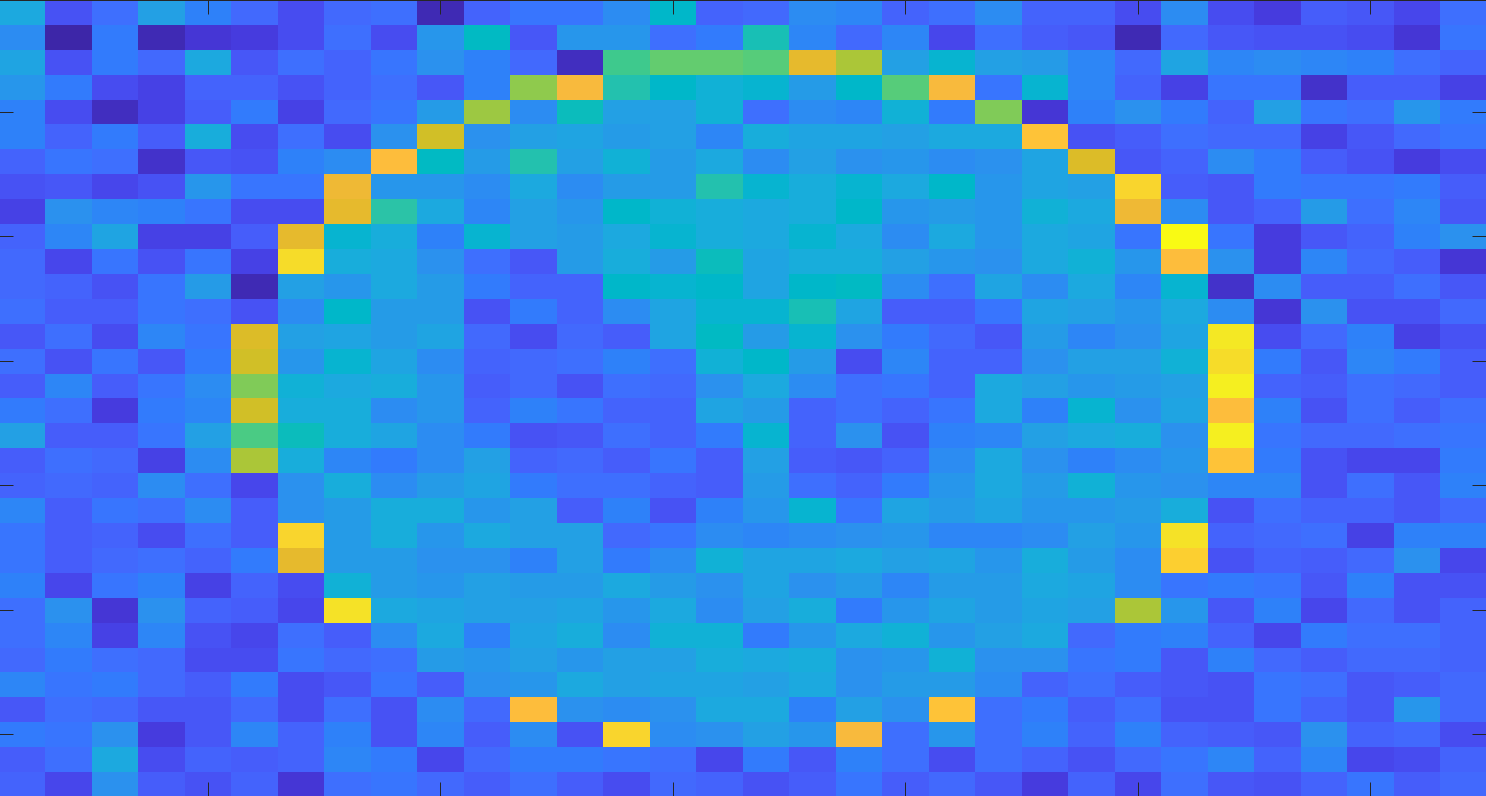}
    &
    \includegraphics[width=0.28\textwidth,height=1.35in]{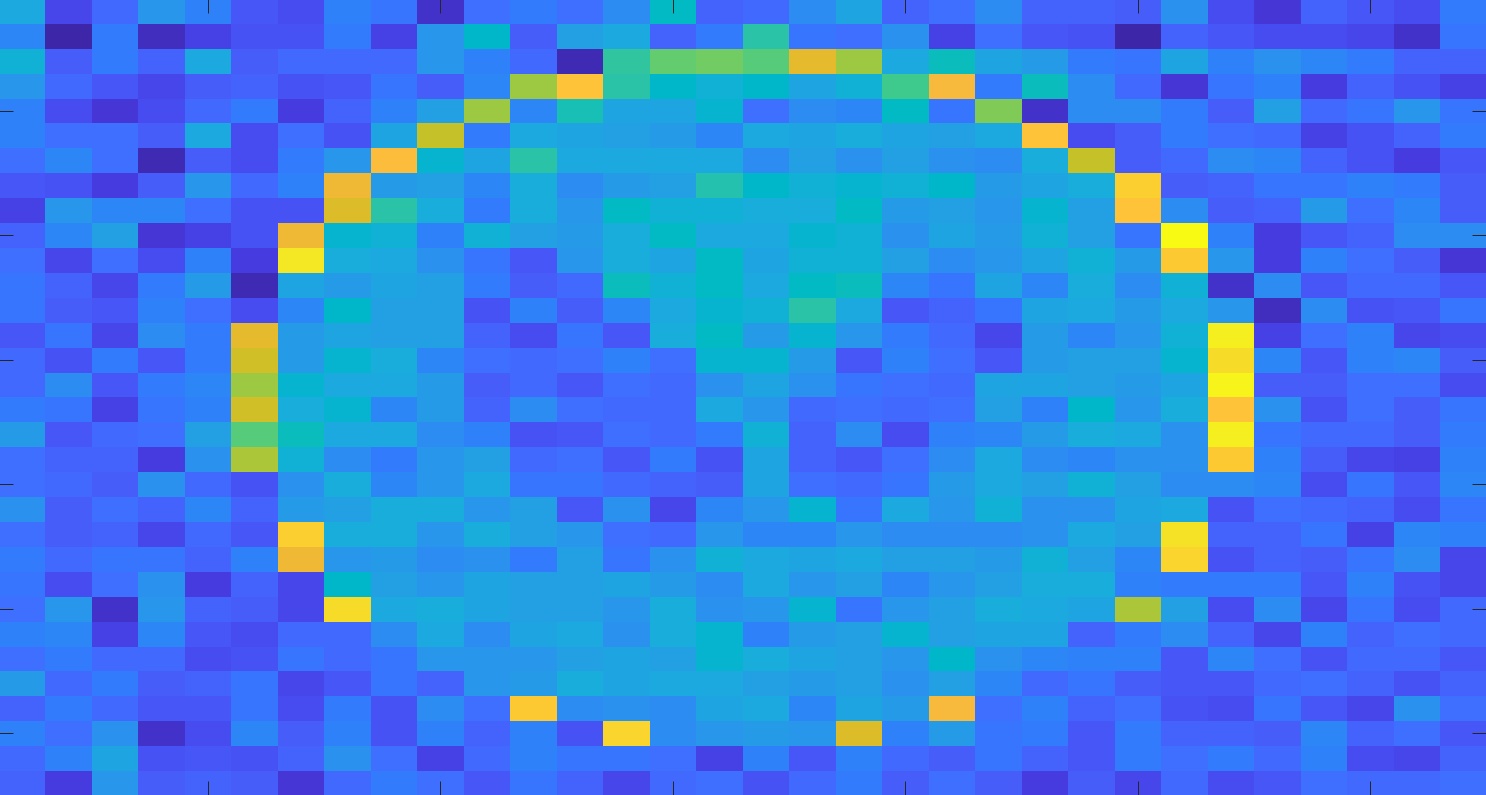}
    \\
    & & \textbf{relerr} $\approx $ 3.49e-01& \textbf{relerr} $\approx$ 3.75e-01
  \end{tabular}
  \begin{tabular}{ccccc}
    & $\bfW_1$ & $\bfW_2$ & $\bfW_3$ & $\bfW_4$ \hspace{-5mm}
    \\
    \rottextt{\small{\textbf{Weights}}} \hspace{-5mm}
    &
    \includegraphics[width=0.205\textwidth,height=1.2in]{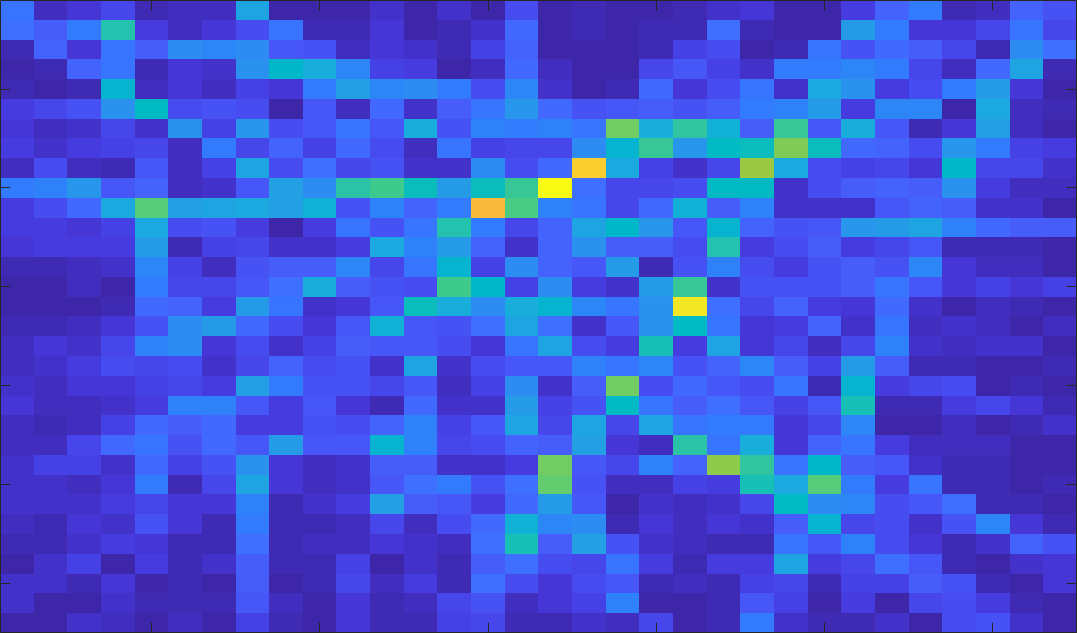}
    &
    \includegraphics[width=0.205\textwidth,height=1.2in]{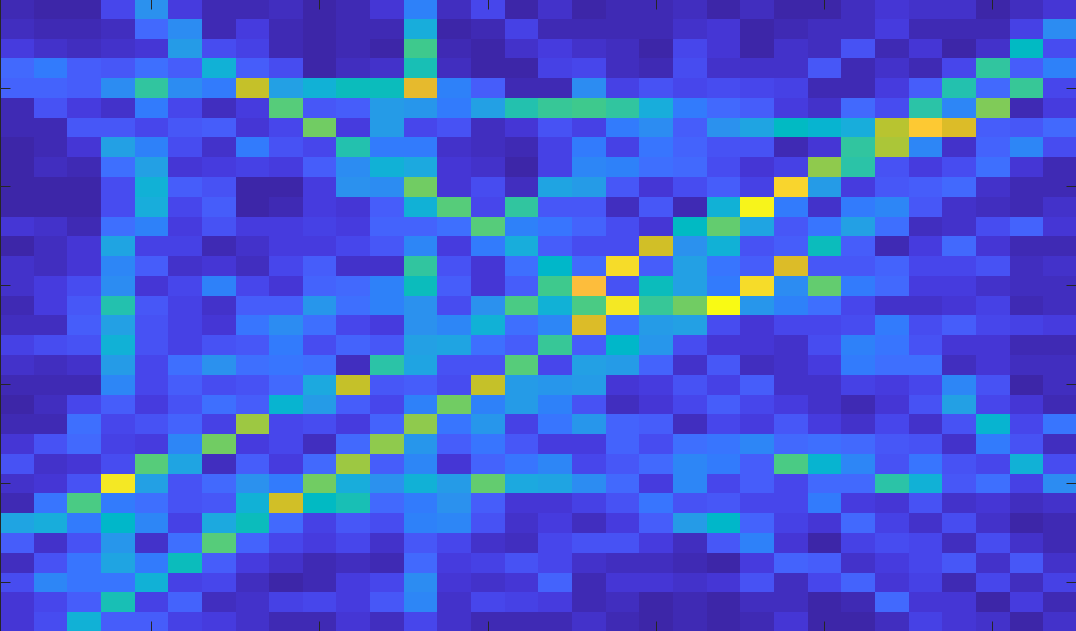}
    &
    \includegraphics[width=0.205\textwidth,height=1.2in]{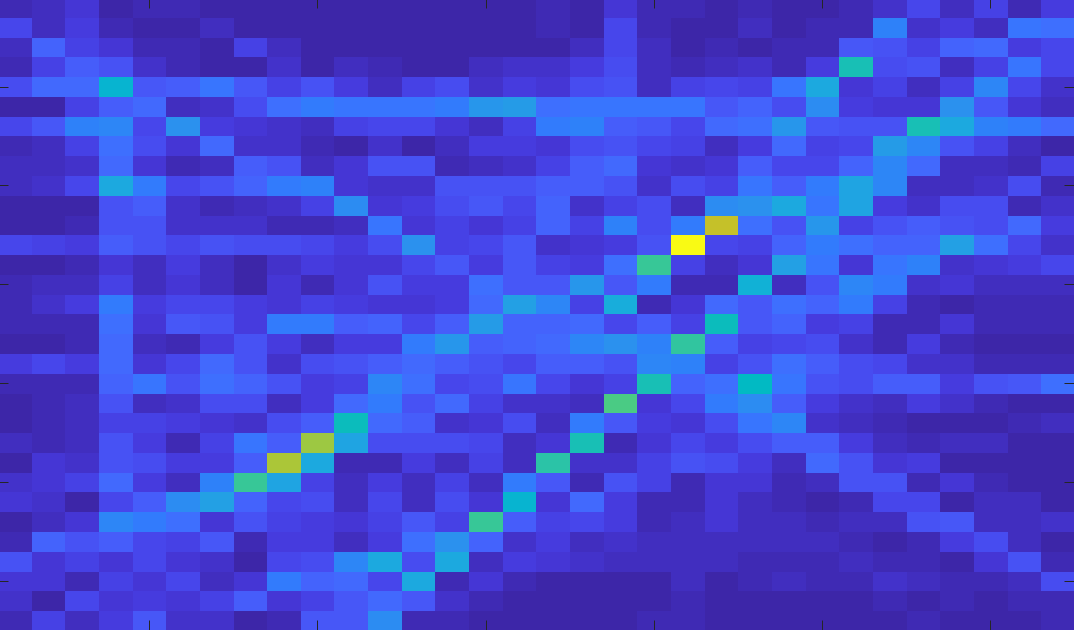}
    &
    \includegraphics[width=0.205\textwidth,height=1.2in]{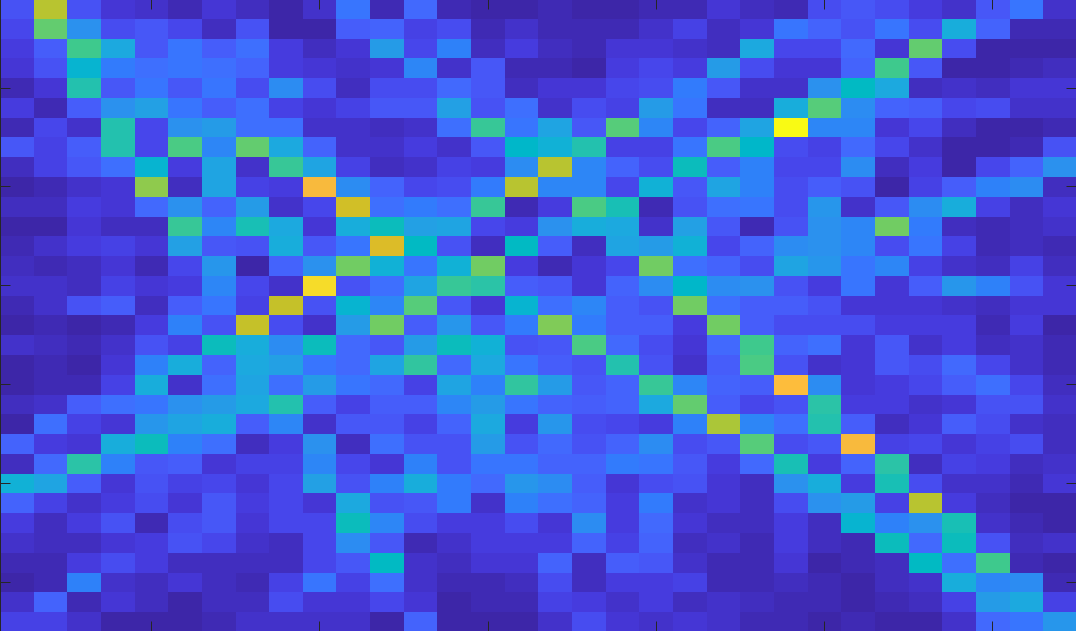}   
  \end{tabular}
  \caption{\textit{Reconstructions (first row) after 10 iterations and weights (second row) for the tomography problem from \texttt{\emph{Regtools~\cite{hansen1994regularization}}}. 
  }}
  \label{fig:tomo}
  \vspace{-5mm}
\end{figure}
\begin{table}[!t]
\centering
\small
\begin{tabular}{|c|c|c|c|c|c|}
\hline
\multicolumn{6}{|c|}{\textbf{UF Sparse Matrix Collection Results}}
\\
\hline
& & \multicolumn{2}{|c|}{\textbf{unweighted ADMM}} & \multicolumn{2}{c|}{\textbf{weighted ADMM}}
\\
\hline
\textbf{matrix}        & \textbf{cond \#}     & \textbf{residual}       & \textbf{rel. error}     & \textbf{residual}      & \textbf{rel. error}   \\ \hline
bcspwr03              & 5.01e+02             & 4.47e-02                 & 2.58e-01                 & 1.97e-02             & 1.63e-01          \\ \hline
bcsstk03              & 6.79e+06             & 6.67e-01                 & 9.99e-01              & 5.82e-01             & 9.91e-01             \\ \hline
bcsstk19              & 1.34e+11             & 2.81e-01                 & 9.01e-01              & 6.30e-02             & 8.62e-01             \\ \hline
bfwb782               & 1.81e+01             & 4.94e-02                 & 1.04e-01              & 2.83e-01             & 1.38e-01             \\ \hline
can\_229              & 4.01e+17             & 5.13e-02                 & 2.10e-01              & 1.81e-02             & 1.46e-01             \\ \hline
cavity02              & 8.12e+04             & 6.07e-01                 & 9.33e-01              & 2.46e-01             & 7.83e-01             \\ \hline
cavity03              & 5.85e+05             & 5.69e-01                 & 9.01e-01              & 1.90e-01             & 7.34e-01             \\ \hline
ch5-5-b4              & 1.00e+00             & 1.21e-01                 & 9.85e-01              & 5.01e-03             & 9.84e-01             \\ \hline
dwt\_307              & 2.35e+18             & 8.92e-02                 & 1.82e-01              & 2.23e-02             & 9.20e-02             \\ \hline
football              & 3.74e+02             & 5.90e-02                 & 4.88e-01              & 2.44e-02             & 3.58e-01             \\ \hline
fs\_183\_3            & 3.27e+13             & 7.33e-02                 & 1.00e+00              & 2.32e-02             & 1.00e+00             \\ \hline
G23                   & 1.00e+04             & 2.27e-02                 & 2.63e-01              & 1.88e-02             & 2.55e-01         \\ \hline
GD98\_c               & 9.87e+16             & 7.78e-02                 & 3.70e-01              & 5.13e-02             & 2.62e-01         \\ \hline
gre\_115              & 4.97e+01             & 2.77e-01                 & 4.70e-01              & 7.64e-02             & 2.67e-01         \\ \hline
gre\_343              & 1.12e+02             & 1.18e-01                 & 1.77e-01              & 4.48e-02             & 6.90e-02         \\ \hline
grid1\_dual           & 3.35e+16             & 3.74e-02                 & 3.67e-01              & 2.56e-02             & 3.20e-01         \\ \hline
impcol\_d             & 2.06e+03             & 3.50e-01                 & 7.07e-01              & 9.74e-02             & 3.94e-01         \\ \hline
jpwh\_991             & 1.42e+02             & 1.89e-01                 & 8.81e-01              & 1.61e-01             & 8.66e-01         \\ \hline
lowThrust\_1          & Inf                  & 4.40e-01                 & 9.96e-01              & 2.84e-01             & 9.82e-01         \\ \hline
lund\_a               & 2.80e+06             & 1.06e-01                 & 5.99e-01              & 4.05e-02             & 5.65e-01         \\ \hline
nos3                  & 3.77e+04             & 5.27e-01                 & 9.88e-01              & 2.22e-01             & 9.67e-01         \\ \hline
odepa400              & 2.26e+05             & 4.91e-01                 & 9.99e-01              & 1.86e-01             & 9.97e-01         \\ \hline
pde900                & 1.53e+02             & 5.98e-01                 & 9.82e-01              & 2.90e-01             & 9.39e-01         \\ \hline
poisson2D             & 1.33e+02             & 3.33e-01                 & 7.44e-01              & 8.03e-02             & 6.64e-01         \\ \hline
polbooks              & 7.20e+02             & 3.72e-02                 & 2.97e-01              & 2.30e-02             & 2.51e-01         \\ \hline
problem1              & 3.11e+16             & 4.08e-01                 & 9.11e-01              & 1.43e-01             & 8.24e-01         \\ \hline
rdb200l               & 1.33e+02             & 8.68e-02                 & 1.44e-01              & 1.91e-02             & 1.05e-01         \\ \hline
str_0                 & 2.74e+02             & 8.46e-01                 & 9.37e-01              & 4.04e-01             & 7.02e-01         \\ \hline
TF10                  & 7.34e+02             & 2.63e-01                 & 4.82e-01              & 5.21e-02             & 2.85e-01         \\ \hline
young1c               & 4.15e+02             & 6.39e-01                 & 9.40e-01                  & 3.59e-01             & 6.97e-01         \\ \hline
\end{tabular}
  \caption{\textit{Comparison of the accuracy obtained using the unweighted and weighted ADMM applied to least-squares problems from the UF Sparse Matrix Collection~\cite{davis2011university}. Columns 1 and 2 show the name and condition number of the matrices. Columns 3 an 4 show the relative residuals of the unweighted and weighted ADMM at iteration 10, respectively. Columns 5 and 6 show the relative errors of the unweighted and weighted ADMM at iteration 10, respectively. 
  }}
  \label{tab:UFLeastSquares}
\end{table}
\begin{table}[!t]  
  \centering
  \begin{tabular}{ccc|ccc}
    \multicolumn{2}{c}{ \textbf{a) DC Resistivity}} & & & \multicolumn{2}{c}{ \textbf{b) Travel-Time Tomography}}
    \\
    \hline
          $\nabla \cdot (\sigma(x_{DC}) \nabla u) = q$
          & in $ \Omega$ 
          &
          &
          & $| \nabla u |^2 = x_{Eik}$
          & in $ \Omega$
    \\
          $\nabla u \cdot \vec{n} = 0$
          & $ \text{on } \partial \Omega$
          &
          &
          & $u(x_0) = 0$
          & 
    \\
          $ u \to 0$
          & $x \to \infty$
          &
          &
          &
          &
   \end{tabular}
   \caption{PDEs corresponding to two different geophysical imaging techniques: 1) DCR (left), where $\sigma$ denotes the ground conductivity parametrized by our model $x_{DC} \in \Omega$, $u: \Omega \mapsto R$ is the electric potential field, $q: \Omega \mapsto \bbR$ are the sources, and 2) travel-time tomography (right) modeled using the Eikonal equation, where $u: \Omega \mapsto \bbR$ is the arrival time of the first wave that evolves from the source $q_0$ located at $x_0$, and $x_{Eik}: \Omega \mapsto \bbR$ is the squared slowness.}
   \label{tab:PDEs}
   \vspace{-3mm}
\end{table} 
We begin by comparing the weighted and unweighted consensus ADMM on a series of linear least squares problems from \texttt{Regtools}~\cite{hansen1994regularization} and the UF Library of Sparse Matrices~\cite{davis2011university}. For these problems, we use $N=4$ splittings and solve
\begin{equation}
  \begin{split}
    &\argmin_{\bfx_j, \bfz} \;\; \sum_{j=1}^4 \left( \hf \| \bfA_j \bfx_j - \bfy_j \|_2^2 + \frac{\alpha}{2} \| \bfx_j \|_2^2 \right)
    \\
    &\;\;\text{  s.t. } \quad \bfW_j(\bfx_j - \bfz) = \mathbf{0}, \;\;\; j=1,\ldots,4,
  \end{split}
\end{equation}
where similar to Ex. \ref{ex:identity}, $\bfA_j \in \bbR^{(m/4) \times n}$ and $\bfy_j \in \bbR^{m/4}, \; j=1,\ldots,4,$ are chosen by partitioning the rows of the original matrix and the data, $\bfA \in \bbR^{m \times n}$ and $\bfy \in \bbR^m$, respectively. For the deblurring and tomography problems from \texttt{Regtools}, we use the same splittings as in Fig.~\ref{fig:identityExampleWeighted} where we split the rows corresponding to the different quadrants of the image. For the non-image based problems from the UF library, $\bfA_1$ and $\bfy_1$ correspond to the first $m/4$ rows of $\bfA$ and $\bfy$, respectively, $\bfA_2$ and $\bfy_2$ correspond to the second $m/4$ rows of $\bfA$ and $\bfy$, respectively, and so on. In the case that the number of rows, $m$, is not divisible by 4, we round accordingly.

We add a smallness regularization term with $\alpha=10^{-2}$ since the splittings $\bfA_j$ in our experiments are underdetermined ($m/4 < n$), leading to rank-deficient coefficient matrices $\bfA_j^\top \bfA_j$ arising from the normal equations. We set the initial penalty parameter to be $\rho^{(0)} = 5$ and use the adaptive scheme described in~\eqref{eq:adaptiveRho}. We run the unweighted and weighted consensus ADMM for ten iterations and show comparisons of the relative residuals and relative errors. To compute the weights, we follow the procedure in Sec.~\ref{subsec:weights} and compute a rank-10 approximation of the Hessian of the misfits using MATLAB's \texttt{eig} function.

For the \texttt{Regtools} imaging problems, we obtain better results using the weighted ADMM as seen in Fig.~\ref{fig:blur} and Fig.~\ref{fig:tomo}. The effect of the weighting scheme is evident in the lower relative errors and the better reconstructions for the deblurring problem (see Fig.~\ref{fig:blur}). This is in part because the weights in the deblurring problem are very different from one another, leading to better averaging in~\eqref{alg:zstepExplicit}. For the tomography problem, the reconstructions we obtain with the weighted ADMM are only marginally better (see Fig.~\ref{fig:tomo}). This is due to the similarity of the weights in the tomography problem as can be seen in Fig.~\ref{fig:tomo}, which leads to averaging reconstructions that are similar to those of the unweighted ADMM. 

Finally, for the UF matrices, we randomly take $30$ matrices with dimensions $100 \leq m,n \leq 1000$ from the library and compare both methods in Tab.~\ref{tab:UFLeastSquares} after ten iterations. We report their condition number, relative residuals, and relative errors. We obtain better results with the weighted ADMM after ten iterations. 
We refrain from solving these problems in parallel since they are small 2D problems and are mainly used as a proof-of-concept.

\begin{figure}[t]
  \begin{tabular}{ccc}
    (a) Eikonal misfit & & (b) Eikonal relative errors
    \\
    \includegraphics[width=0.45\textwidth, height=1.5in]{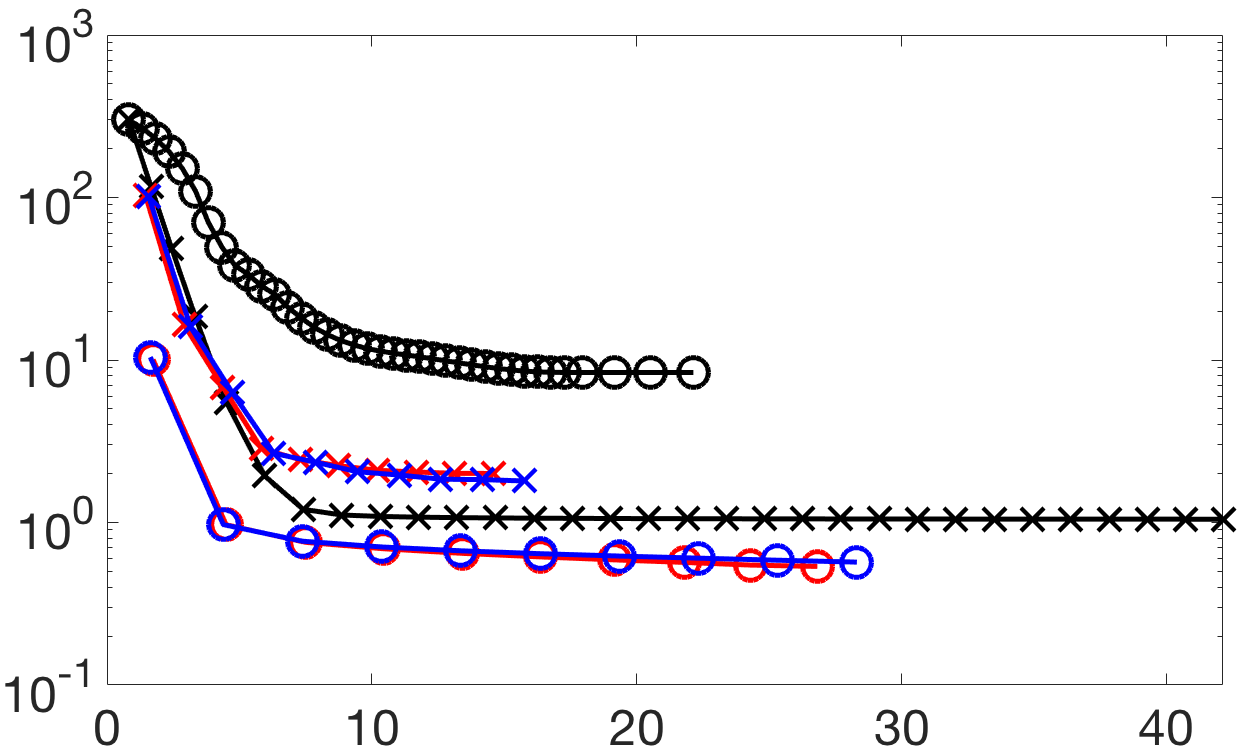}
    &
    &
    \includegraphics[width=0.45\textwidth, height=1.5in]{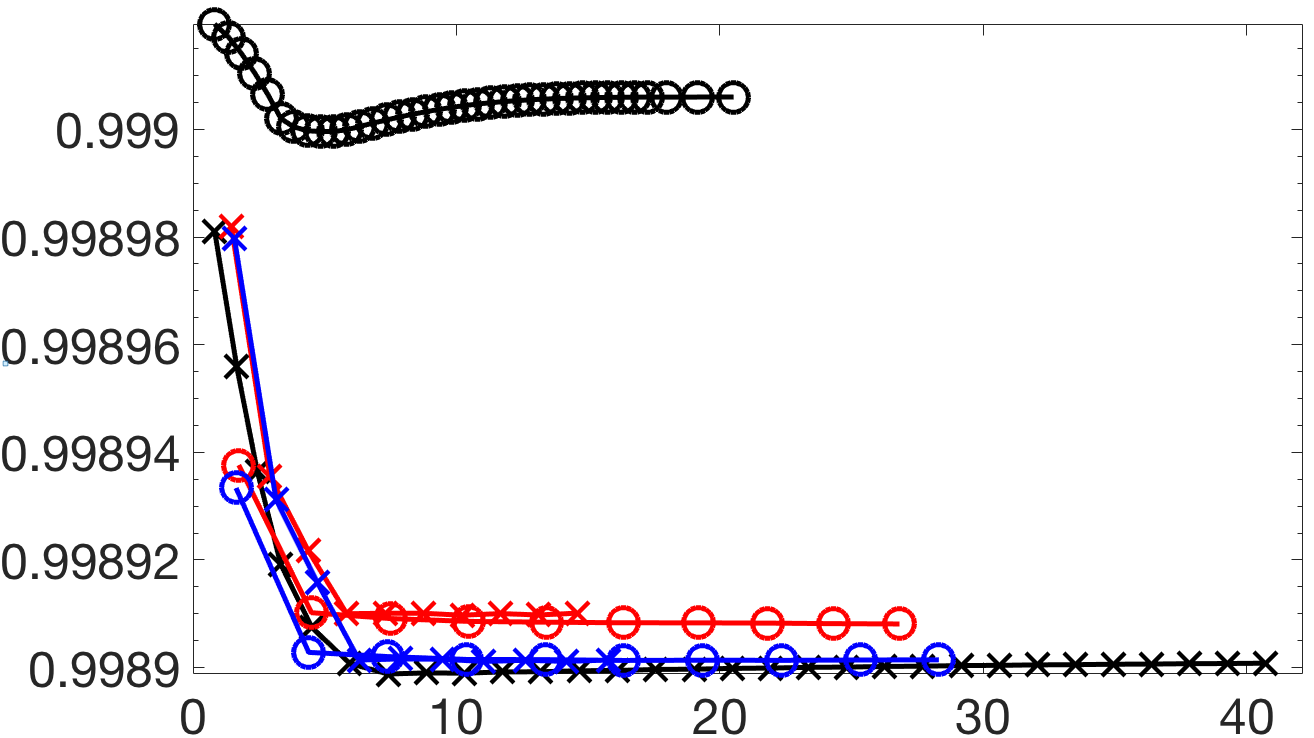}
    \\
    runtime (mins)
    &
    &
    runtime (mins)
    \\
    \\
    (c) joint misfits & & (d) joint relative errors
    \\
    \includegraphics[width=0.45\textwidth, height=1.5in]{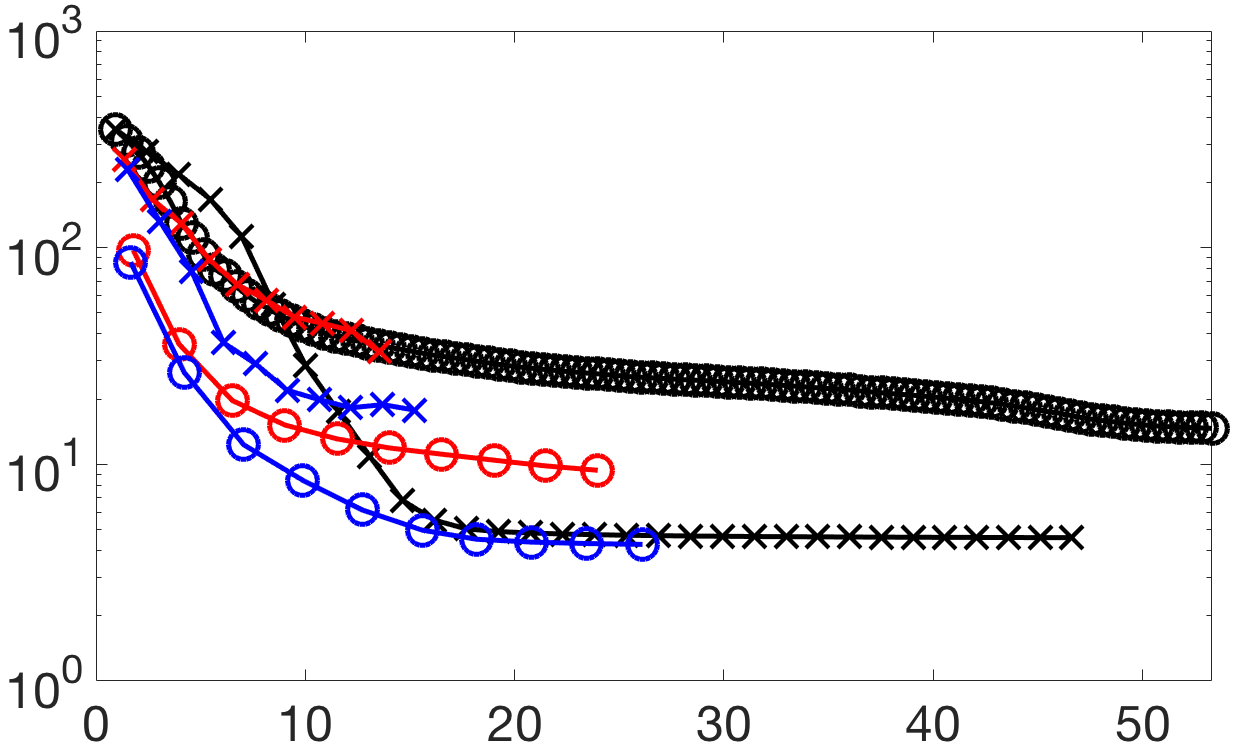}
    &
    &
    \includegraphics[width=0.45\textwidth, height=1.5in]{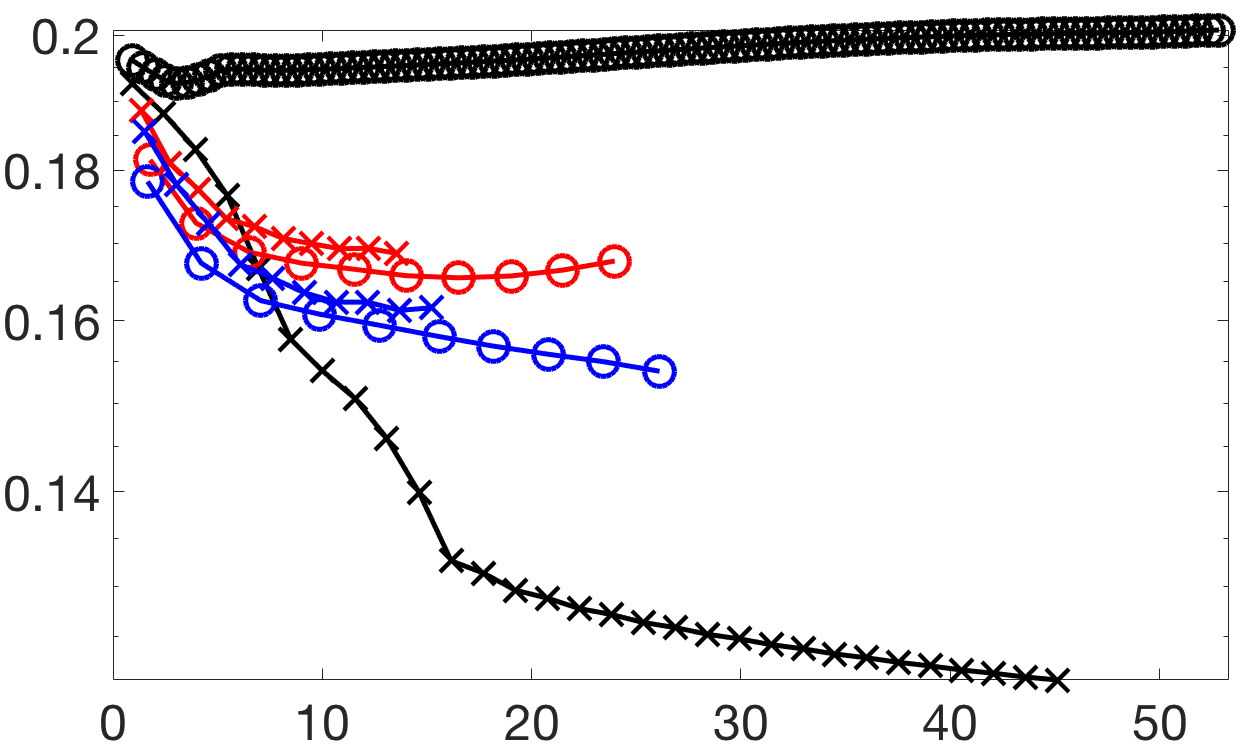}
    \\
    runtime (mins)
    &
    &
    runtime (mins)
    \\
    \\
    \multicolumn{3}{c}{\includegraphics[width=0.6\textwidth]{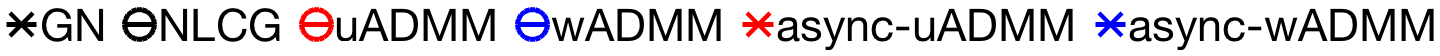}}
  \end{tabular}
  \caption{\textit{Misfit and relative errors for the Eikonal problem using 6 different algorithms: Gauss-Newton, NLCG, uADMM, wADMM, async-uADMM, and async-wADMM. Here, the x-axis represents runtime in minutes. The experiments were run on a shared memory computer operating Ubuntu 14.04 with 2 Intel Xeon E5-2670 v3 2.3 GHz CPUs using 12 cores each, and a total of 128 GB of RAM. Here, Julia is installed and compiled using Intel Math Kernel Library.}}
  \label{fig:MisfitsRelErrs}
  \vspace{-7mm}
\end{figure}
\begin{figure}[!t]
\hspace{-5mm}
\centering
  \begin{tabular}{ccc}
    (a) reference model & (b) ground truth & 
    \\
    \includegraphics[width=0.3\textwidth, height=1.05in]{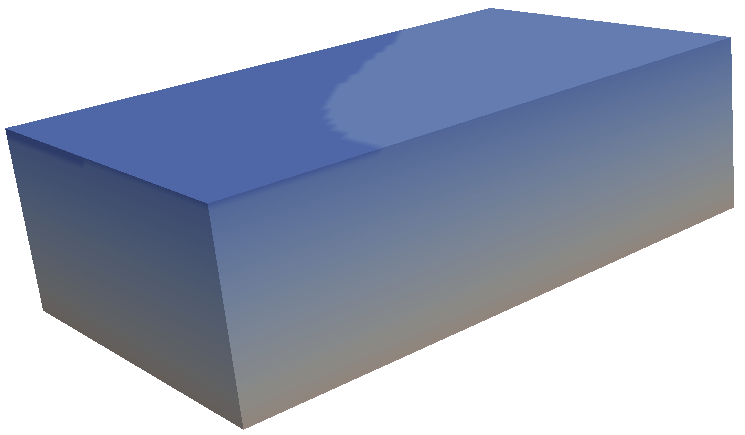}
    &
    \includegraphics[width=0.3\textwidth, height=1.05in]{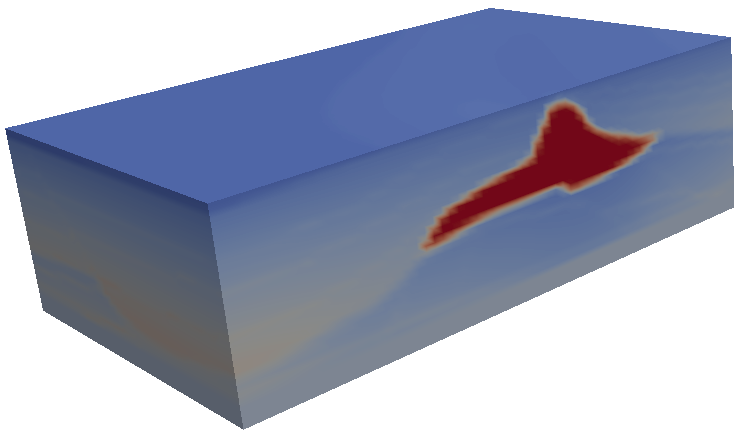}
    &
    \\
    \\
    \\
    & \textbf{Eikonal Reconstructions} &
    \\
    \hline
    \\
    (c) Gauss-Newton & (d) wADMM & (e) uADMM
    \\
    \includegraphics[width=0.301\textwidth, height=1.05in]{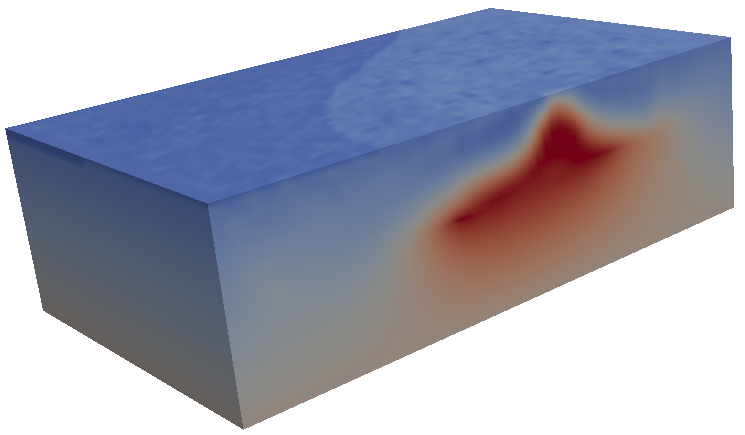}
    &
    \includegraphics[width=0.301\textwidth, height=1.05in]{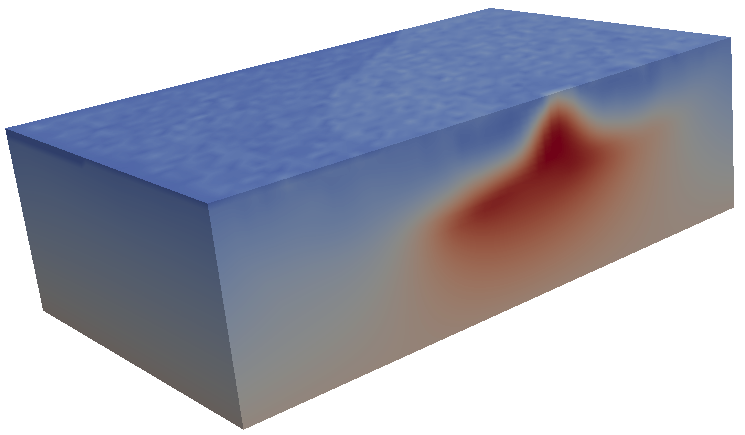}
    &
    \includegraphics[width=0.301\textwidth, height=1.05in]{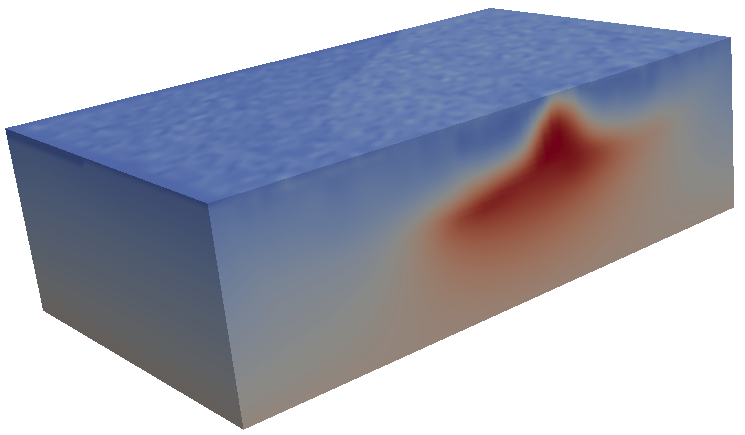}
    \\
    (f) NLCG & (g) async-wADMM & (h) async-uADMM
    \\
    \includegraphics[width=0.301\textwidth, height=1.05in]{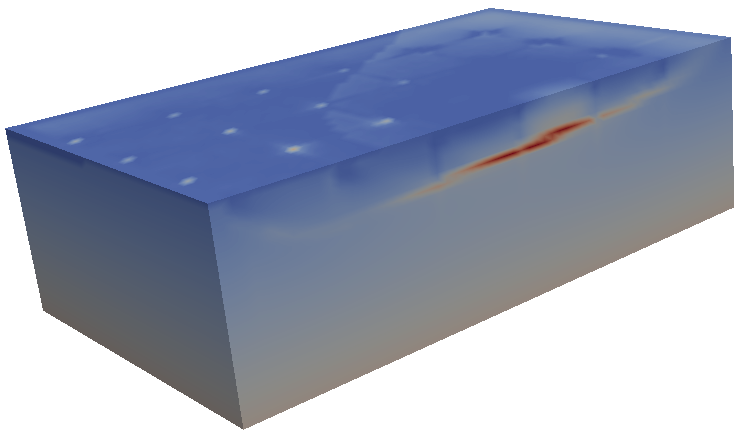}
    &
    \includegraphics[width=0.301\textwidth, height=1.05in]{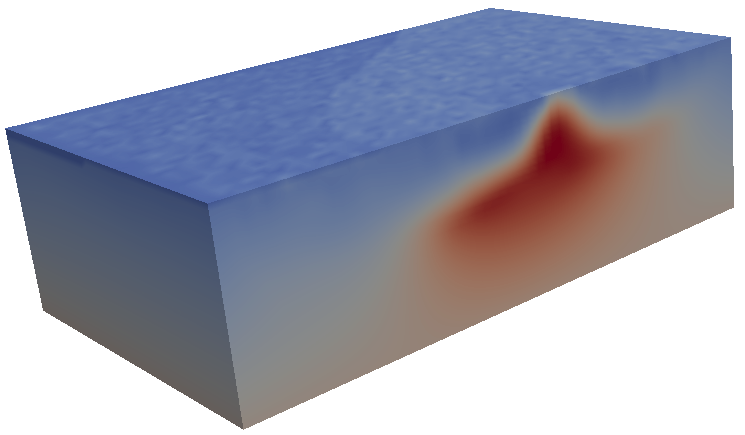}
    &
    \includegraphics[width=0.301\textwidth, height=1.05in]{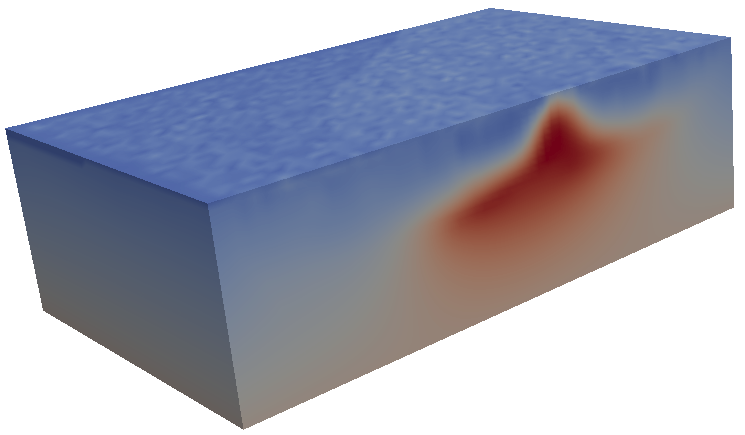}
    \\
    \\
    \\
    & \textbf{Joint Reconstructions} &
    \\
    \hline
    \\
    (i) Gauss-Newton & (j) wADMM & (k) uADMM
    \\
    \includegraphics[width=0.3\textwidth, height=1.05in]{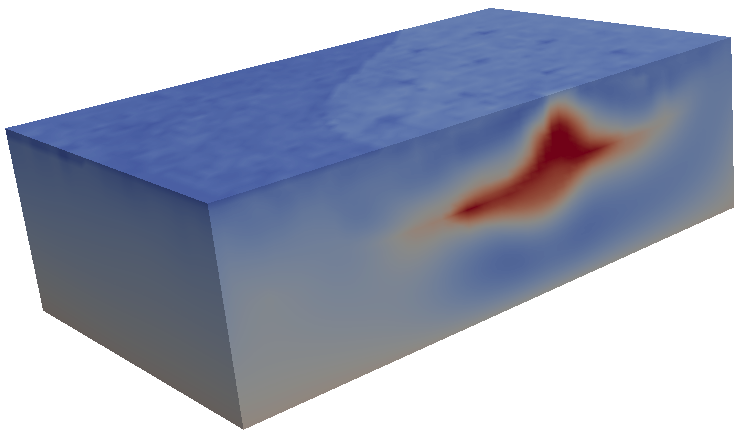}
    &
    \includegraphics[width=0.3\textwidth, height=1.05in]{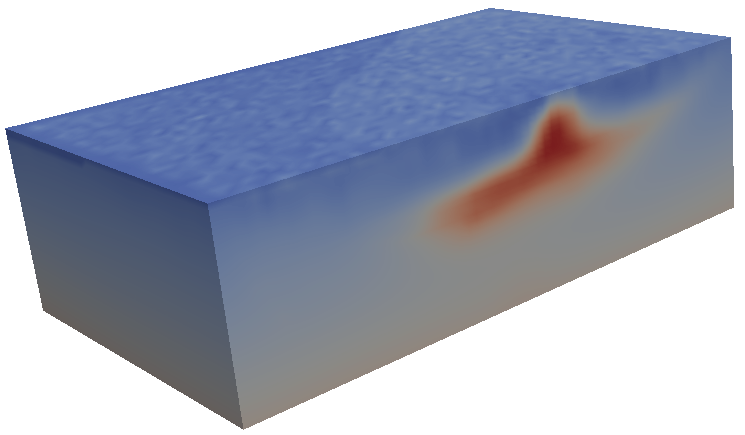}
    &
    \includegraphics[width=0.3\textwidth, height=1.05in]{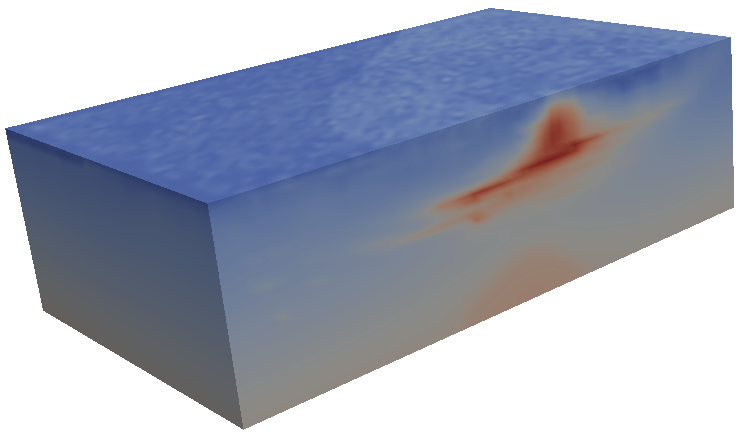}
    \\
    (l) NLCG & (m) async-wADMM & (n) async-uADMM
    \\
    \includegraphics[width=0.3\textwidth, height=1.05in]{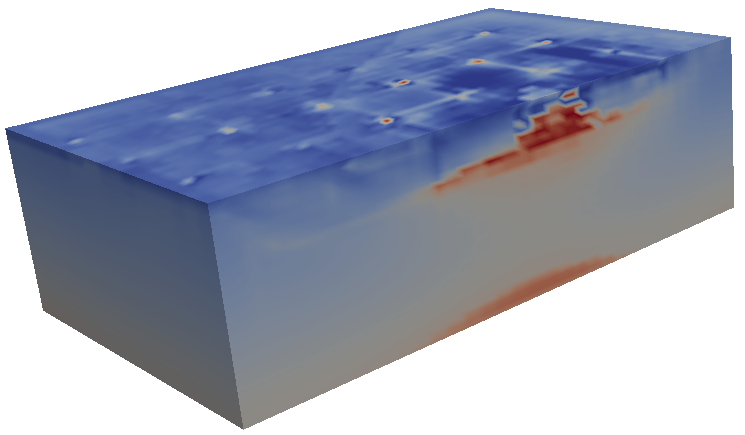}
    &
    \includegraphics[width=0.3\textwidth, height=1.05in]{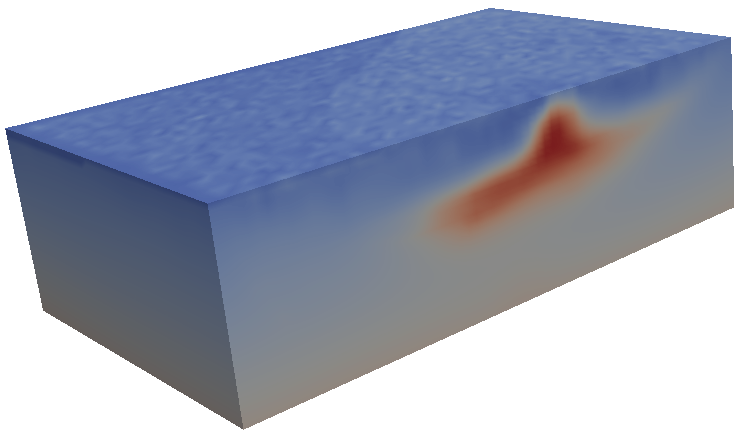}
    &
    \includegraphics[width=0.3\textwidth, height=1.05in]{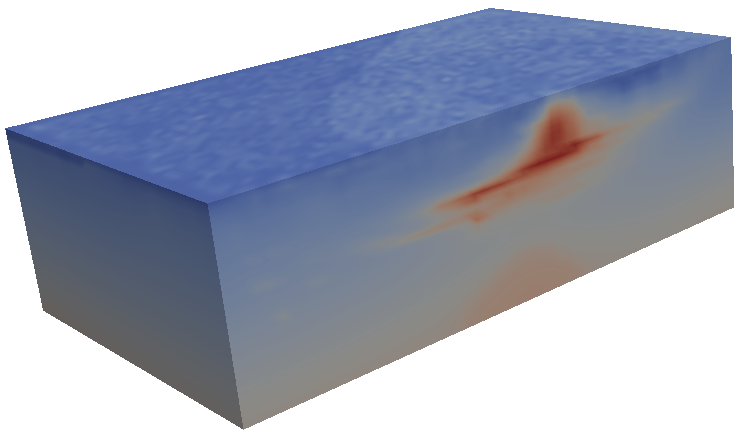}
  \end{tabular}
  \caption{\textit{reconstructions of SEG model with a single-physics and multiphysics experiment.}}
  \label{fig:reconstructions}
  \vspace{-5mm}
\end{figure}
\subsection{Single-Physics Parameter Estimation}
\label{subsec:singlePhysics}
As more a more realistic test problem, we consider the 3D SEG/EAGE model \cite{aminzadeh1997} as the ground truth (see Fig.~\ref{fig:reconstructions}b) and test our method for a single-physics inversion involving the travel-time tomography survey. The model contains a salt dome in which the velocity is significantly higher than in the background. The domain is of size
$13.5$ km $\times$ $13.5$ km $\times$ $4.2$ km and is divided into $64 \times 64 \times 32$ equally sized mesh cells of approximate size of $211 \text{m} \times 211 \text{m}\times 11 \text{m}$ each. We implement our experiments in extension of \texttt{jInv}~\cite{ruthotto2017jinv}, an open-source package for PDE parameter estimation written in Julia \cite{bezanson2017julia}. For brevity, since the travel time tomography problem is modeled by the Eikonal equation, we refer to it as the Eikonal problem for the remainder of the paper. We solve these problems in parallel and experiment on the effect of the asynchronous variant (async-ADMM) on the weighted and unweighted consensus ADMM.

The PDE involved in the forward problem is the Eikonal equation (see Tab.~\ref{tab:PDEs}), and it is solved using the Factored Eikonal Fast Marching Algorithm \cite{treister2016fast}. We solve the inversion using $36$ sources and $3600$ receivers located on the top surface of the domain. We compare the weighted and unweighted ADMM (wADMM and uADMM), their asynchronous variants (async-wADMM and async-uADMM), Gauss-Newton (GN), and NLCG. For all 6 algorithms, we use diffusion regularization with regularization parameter $\alpha = 10^{-3}$ to enforce smoothness. We solve all inversions in parallel using 10 workers. Here, $6$ workers solve forward problems containing $4$ sources each, and the remaining $4$ workers solve forward problems containing $3$ sources each.

We run the Gauss-Newton inversion for a maximum of $30$ outer iterations and use at most 10 PCG iterations with PCG stopping tolerance of $10^{-1}$ to solve the Gauss-Newton system. For the NLCG inversion, we set a maximum of $100$ outer iterations since it is expected to take more iterations than Gauss-Newton to reach the same accuracy. In the ADMM inversions, we run a total of $10$ outer iterations with $3$ GN iterations used to solve the subproblems. This particular choice of inner GN and outer ADMM iterations aims to balance the runtime and computations performed with those of the Gauss-Newton inversion while avoiding solving the subproblems too inexactly (as this may lead to lack of convergence). In the ADMM subproblems, each GN iteration also uses at most $10$ PCG iterations with PCG stopping tolerance of $10^{-1}$ as in the Gauss-Newton inversion. 

For the penalty parameter, we use the scheme described in \eqref{eq:adaptiveRho} to vary $\rho$ and use a lower bound of $10^{-12}$. As expected, the performance of ADMM depends crucially on the initial choice of $\rho$; therefore, we report the best results obtained from initial values of $\rho^{(0)} \in [10^{-8}, 10^2]$. In our experiment, the optimal initial values are $\rho^{(0)} = 10^{-8}$ for uADMM and $\rho^{(0)} = 10^{-2}$ for wADMM. In the asynchronous case, we perform a global update whenever $N_a=5$ workers report their solutions and enforce the bounded delay condition by requiring all workers to report results at least once every $k_a=4$ iterations. To compute the weights, we use the Lanczos bidiagonalization algorithm from \texttt{KrylovMethods} \cite{LarsKrylovMethods} to compute a rank-5 approximation of the approximate Hessians of the data misfits. Again, we note that highly accurate uncertainties are not necessary in our case and a good guess is sufficient for our experiments. The computation of the weights took about 38 seconds.

In Fig.~\ref{fig:MisfitsRelErrs}(a-b), we show the relative errors and misfits for the Eikonal problem. Both wADMM and uADMM outperform the rest of the algorithms in terms of the misfits. However, all algorithms except for NLCG have roughly the same relative errors; this is reflected in the similar reconstructions in Fig.~\ref{fig:reconstructions}(c-h). The impact of communication and latency in the difference of runtimes between asynchronous ADMM variants, which ran for about 15 minutes, and the Gauss-Newton-PCG, which ran for about 43 minutes, is evident.
For the NLCG, a total of $37$ iterations were performed before a linesearch fail was reached. As expected, an iteration from the NLCG method is much quicker than an iteration from the remaining 5 methods since each NLCG iteration only requires explicit steps to update the model.
\subsection{Multiphysics Parameter Estimation}
\label{subsec:jointInversion}
We now add a second modality to Sec.~\ref{subsec:singlePhysics}, the DCR survey, which is modeled by the steady-state heterogeneous diffusion equation (see Tab.~\ref{tab:PDEs}), and consider a multiphysics inversion. Here, we keep the same settings for the Eikonal problem and use $32$ sources and $1682$ receivers located on the top surface of the domain for the DCR survey. To solve the DCR forward problem, we Julia's direct solver. We assume known petrophysics \cite{schon2015physical} so that we have a relation between the ground conductivity $\sigma$ and the wave velocity $x$ given by
\begin{equation}
  \sigma(x) = \left(2 - \frac{x}{c} \right) \left(\frac{b-a}{2} (\tanh(10(c-x))+1)+a \right).
\end{equation}
Here, a and b are the conductivity values set to $0.1$ and $1.0$ respectively, and $c = 3.0$ is the velocity in which the contrast is centered. More details can be found in \cite{ruthotto2017jinv}. 

As in Sec.~\ref{subsec:singlePhysics}, we compare six algorithms: wADMM, uADMM, async-wADMM, async-uADMM, Gauss-Newton-PCG, and NLCG.
We solve all the inversions in parallel using ten workers. We assign the DCR problem to one worker since it is easier to solve and assign the Eikonal problem to the remaining nine workers. In this case, the nine workers in charge of the Eikonal problem solve forward problems containing $4$ sources each. The inversion settings are also the same as in Sec.~\ref{subsec:singlePhysics} except for the choice of initial penalty parameter, where we find the optimal initial values to be $\rho^{(0)} = 10^{-6}$ for uADMM and $\rho^{(0)} = 1.0$ for wADMM. We also follow the same procedure as in Sec.~\ref{subsec:singlePhysics} to compute the weights for this setup, which took about $54$ seconds.

The results for the relative errors and misfits for the joint inversion can be seen in Fig.~\ref{fig:MisfitsRelErrs}(c-d). Here, the weighted scheme gives us a substantial improvement in our relative error and misfits; this is also reflected in the reconstructions shown in Fig.~\ref{fig:reconstructions}(i-n), where the wADMM and async-wADMM give more accurate reconstructions. We also save on latency and communication with the async-wADMM, which ran for about (16 minutes) while maintaining a good quality of the reconstruction. As expected, the joint inversions enhance the quality of the reconstruction since the different physics involved capture different properties of the model \cite{ruthotto2017jinv}.

\section{Conclusion} \label{sec:conclusion} 
We propose a weighted asynchronous consensus ADMM (async-wADMM) method for solving large-scale PDE parameter estimation problems in parallel.
To this end, the data involved in the problem is divided among the available workers. 
Our scheme is geared toward applications such as PDE parameter estimation where only a few iterations can be afforded. 
Our proposed weighting scheme improves the convergence of the standard ADMM. 
Since our weights are informed by the uncertainties in the estimates of the subproblems in \eqref{alg:xstepExplicit}, we formulate the parameter estimation problem in a Bayesian setting. It is important to note that our scheme can also be applied in the frequentist setting as long as weights are available.  To obtain an overall efficient scheme, we follow the works of \cite{flath2011fast} to quantify the uncertainties in a tractable manner.

As test problems, we solve a collection of linear least-squares problems for proof-of-concept as well as a more realistic single-physics involving the travel time tomography survey, and a multiphysics parameter estimation problem involving the DCR and travel time tomography survey.
Our numerical results show that our method accelerates the convergence of consensus ADMM, particularly in the early iterations.
The quality of the parameter estimate obtained by the weighted async-ADMM scheme is comparable to those of the Gauss-Newton-PCG method; however, the weighted async-ADMM method requires substantially less communication among workers and has smaller latencies, resulting in reduced inversion runtimes. Moreover, since we can choose any optimization scheme to solve the subproblems in async-ADMM, the method sits at a higher level of abstraction and provides additional flexibility. Each subproblem can, therefore, be solved with a tailored solver, making the weighted async-ADMM especially attractive for large-scale multiphysics PDE parameter estimation problems.
For brevity, we do not show the case where the weights are computed in every iteration; however, in this case, we obtain indiscernible reconstructions from those shown in Fig.~\ref{fig:reconstructions}. We intend to further explore our method for large-scale problems where the Gauss-Newton-PCG method cannot be used as well as on computational environments with small communication bandwidth such as cloud computing platforms.

\section*{Acknowledgments}
This material is supported by the U.S. National Science Foundation (NSF) through awards DMS 1522599 and DMS 1751636

\FloatBarrier
\bibliographystyle{abbrv}
\bibliography{references.bib}

\end{document}